\documentclass[a4paper,12pt,reqno]{amsart}

\usepackage[english]{babel}
\usepackage[utf8]{inputenc}
\usepackage[T1]{fontenc}

\usepackage[a4paper,top=2cm,bottom=2cm,left=2cm,right=2cm]{geometry}

\usepackage{titlesec}
\usepackage{lipsum}
\usepackage{enumerate}
\usepackage{enumitem}
\usepackage{mathrsfs}
\usepackage{amsmath,amssymb,amsthm}
\usepackage{mathtools}
\usepackage{graphicx,float,subfigure}
\usepackage{url}
\usepackage{bbold}
\usepackage{color}
\usepackage{xcolor}

\usepackage[all,cmtip]{xy}
\usepackage{multirow}
\usepackage[colorinlistoftodos]{todonotes}
\usepackage[colorlinks=true, allcolors=black]{hyperref}
\usepackage{harpoon}
\usepackage{MnSymbol}

\DeclareMathOperator{\sgn}{sgn}
\DeclareMathOperator{\wt}{wt}
\numberwithin{equation}{subsection}

\theoremstyle{plain}
\newtheorem{theorem}{\scshape Theorem}[section]

\newtheorem{lemma}[theorem]{\scshape Lemma}
\newtheorem{corollary}[theorem]{\scshape Corollary}

\newtheorem*{assumption*}{\scshape Assumption}

\newtheorem*{claim*}{Claim}

\theoremstyle{definition}

\newtheorem{definition}[theorem]{\scshape Definition}
\newtheorem{remark}[theorem]{\scshape Remark}
\newtheorem{example}[theorem]{\scshape Example}

\newcommand{\M}{\operatorname{M}}

\counterwithin{equation}{section}

\titleformat{\section}{\centering\bfseries}{\thesection}{1em}{\MakeUppercase}
\begin{document}
\title{An extension of the Lindstr{\"o}m-Gessel-Viennot theorem}
\author{Yi-Lin Lee}
\address{Department of Mathematics, Indiana University, Bloomington, Indiana 47405}
\email{yillee@iu.edu}

\maketitle
\begin{abstract}
  Consider a weighted directed acyclic graph $G$ having an upward planar drawing. We give a formula for the total weight of the families of non-intersecting paths on $G$ with any given starting and ending points. While the Lindstr{\"o}m-Gessel-Viennot theorem gives the signed enumeration of these weights (according to the connection type), our result provides the straight count, expressing it as a determinant whose entries are signed counts of lattice paths with given starting and ending points.
\end{abstract}

\section{Introduction}

In algebraic and enumerative combinatorics, the Lindstr{\"o}m-Gessel-Viennot theorem is a powerful and elegant result with numerous applications in different contexts. This result was independently discovered by Gessel and Viennot~\cite{GV85}, by Lindstr{\"o}m~\cite{Lin73}, and even earlier in the work on coincidence probabilities by Karlin and McGregor~\cite{KM59}.

Throughout this paper, by a directed graph $G=(V,E)$ (or digraph) we mean a locally finite, simple, connected graph, with a weight function $\wt: E \rightarrow \mathcal{R}$ that assigns elements in some commutative ring $\mathcal{R}$ to each edge of $G$. For brevity, we sometimes call them simply digraphs. A digraph is \textit{acyclic} if it has no directed cycles.

The \textit{weight} of a path $p$ is the product $\wt(p) = \prod_{e}\wt(e)$, where the product is over all edges $e$ of the path $p$. A path with length zero has weight $1$ by convention. The weight of an $n$-tuple of paths $P=(p_1,\dotsc,p_n)$ is the product of the weights of each path: $\wt(P) = \prod_{i=1}^{n} \wt(p_i)$. Let $\mathscr{P}$ be a collection of (families of) paths, we will write $GF(\mathscr{P})$ for the generating function according to the weight $\wt$, that is, $GF(\mathscr{P}) = \sum_{P \in \mathscr{P}} \wt(P)$. In particular, we define
\begin{equation}\label{gf}
  h(u,v) = GF(\mathscr{P}(u,v)) = \sum_{p \in \mathscr{P}(u,v)}\wt(p),
\end{equation}
where $\mathscr{P}(u,v)$ is the set of all paths connecting the vertex $u$ to the vertex $v$.

The general version of the Lindstr{\"o}m-Gessel-Viennot theorem on a directed acyclic graph $G$ is stated below in Theorem~\ref{thm.GV1}, while the simplified but most commonly used version is given in Corollary~\ref{cor.GV2}.
\begin{theorem}[Lindstr\"om~\cite{Lin73}; Gessel,Viennot~\cite{GV85}; Stembridge~\cite{Stem90}]\label{thm.GV1}
   Suppose $U = \{u_1,\dotsc,u_n\}$ and $V = \{v_1,\dotsc,v_n\}$ are two distinct sets of vertices of $G$. Let $\mathscr{P}^{\pi}_{0}(U,V)$ be the set of $n$-tuples of non-intersecting paths $(p_1,\dotsc,p_n)$ with the connection type $\pi \in \mathfrak{S}_n$, namely, the path $p_i$ goes from $u_i$ to $v_{\pi(i)}$, for $1 \leq i \leq n$. Then
  \begin{equation}\label{eq.GV}
    \sum_{\pi \in \mathfrak{S}_n}\sgn(\pi) GF(\mathscr{P}^{\pi}_{0}(U,V)) = \det \left( h(u_i,v_j) \right)_{i,j=1}^{n}.
  \end{equation}
\end{theorem}
We say two sets of the vertices $U = \{u_1,\dotsc,u_n\}$ and $V = \{v_1,\dotsc,v_n\}$ of $G$ are \textit{compatible} if any pair of paths $p_i$ from $u_i$ to $v_k$ and $p_j$ from $u_j$ to $v_{\ell}$ with $i<j,k>\ell$, the paths $p_i$ and $p_j$ must intersect. The key point for this condition is that when $U$ and $V$ are compatible, the $n$-tuples of non-intersecting paths only consist of paths connecting $u_i$ to $v_i$ for $i=1,\dotsc,n$. This is the case when $\pi = \mathsf{id}$ in Theorem~\ref{thm.GV1} which leads to the following corollary.
\begin{corollary}\label{cor.GV2}
  If $U = \{u_1,\dotsc,u_n\}$ and $V = \{v_1,\dotsc,v_n\}$ are compatible, then we have
  \begin{equation}\label{eq.cor}
    GF(\mathscr{P}^{\mathsf{id}}_{0}(U,V)) = \det \left( h(u_i,v_j) \right)_{i,j=1}^{n}.
  \end{equation}
\end{corollary}
In~\cite{Stem90}, Stembridge generalized this result to the situation when the endpoints of the paths are not fixed, and this leads to Pfaffian expressions. We refer the reader to the survey paper written by Krattenthaler~\cite[Section 10.13]{Kra17} for more details.

The Lindstr{\"o}m-Gessel-Viennot theorem is relevant in the enumeration of the semi-standard Young tableaux and their variations, the enumeration of various types of plane partitions, and in evaluating some special kinds of Hankel determinants (see for example~\cite{GV85,GV89},~\cite{Stem90},~\cite{Aig01}). A combinatorial proof of the Jacobi-Trudi type identities for Schur functions is another application (see~\cite[Chapter 4]{Sagan}). In particular, it plays an important role in the enumeration of tilings.

A \textit{tiling} is a covering of a given region on the plane using a given set of tiles without gaps or overlaps. Finding formulas that count (weighted) tilings of specific regions is an important but difficult problem in the field of enumerative combinatorics. The path method (see~\cite[Section 3.1]{Propp15}) is one of the powerful techniques used to count (weighted) tilings. The core idea is to view a tiling as a family of non-intersecting paths and then use the Lindstr{\"o}m-Gessel-Viennot theorem to obtain the number of these families by evaluating a determinant.

For example, the number of domino tilings of the Aztec diamond of order $n$ is given by the elegant expression $2^{n(n+1)/2}$, first proved in~\cite{ELKP1,ELKP2}. Bijective proofs using the path method appeared in~\cite{BK05},~\cite{EF05} and later in~\cite{BvL13}. Domino tilings are in one-to-one correspondence with families of non-intersecting Schr{\"o}der paths; evaluating the determinant of a matrix whose entries are large Schr{\"o}der numbers gives then this elegant formula.

\begin{figure}[hbt!]
\centering
\subfigure[]
{\label{Mbox}\includegraphics[height=0.32\textwidth]{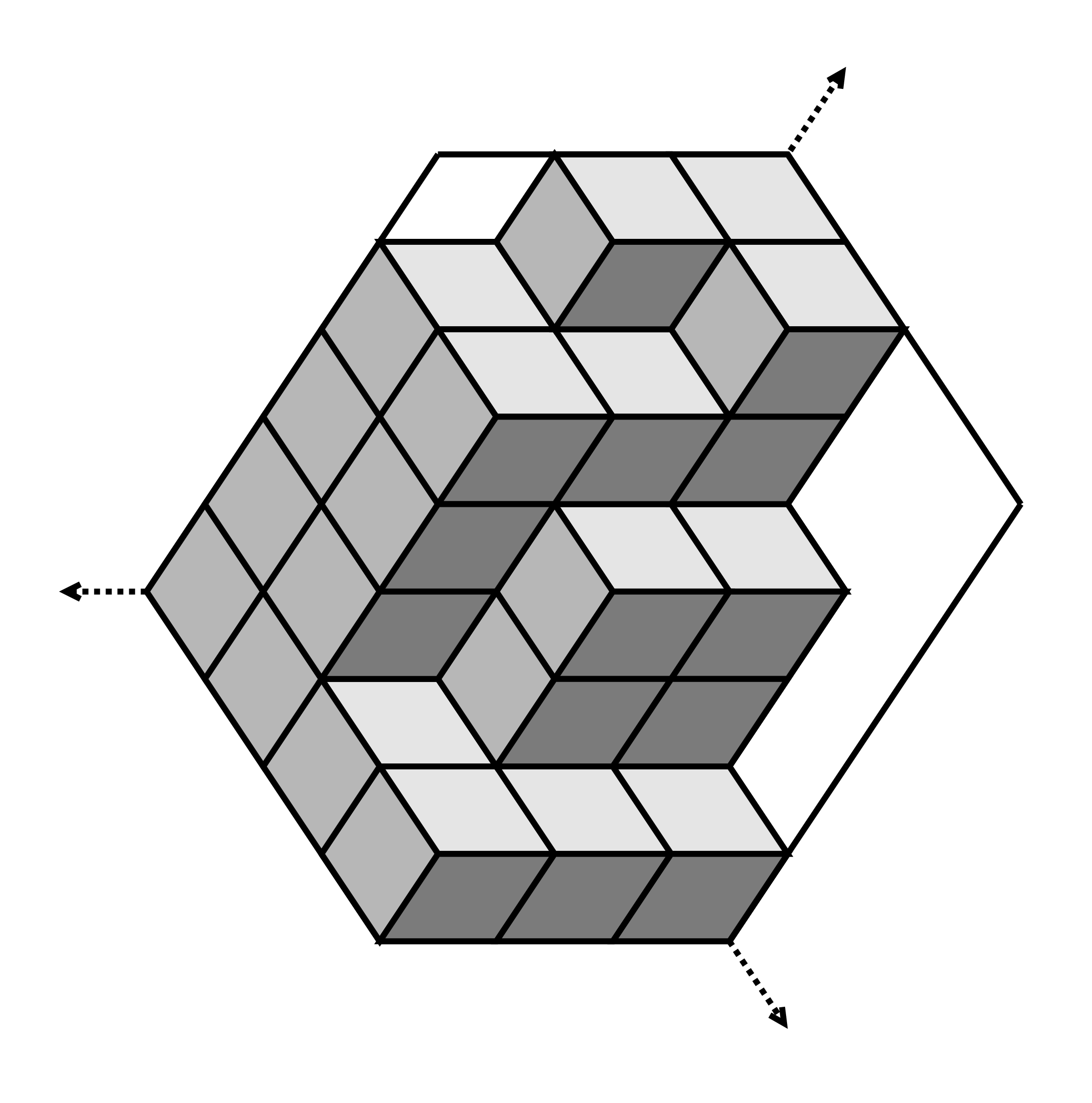}}
\hspace{1mm}
\subfigure[]
{\label{Mtile}\includegraphics[height=0.29\textwidth]{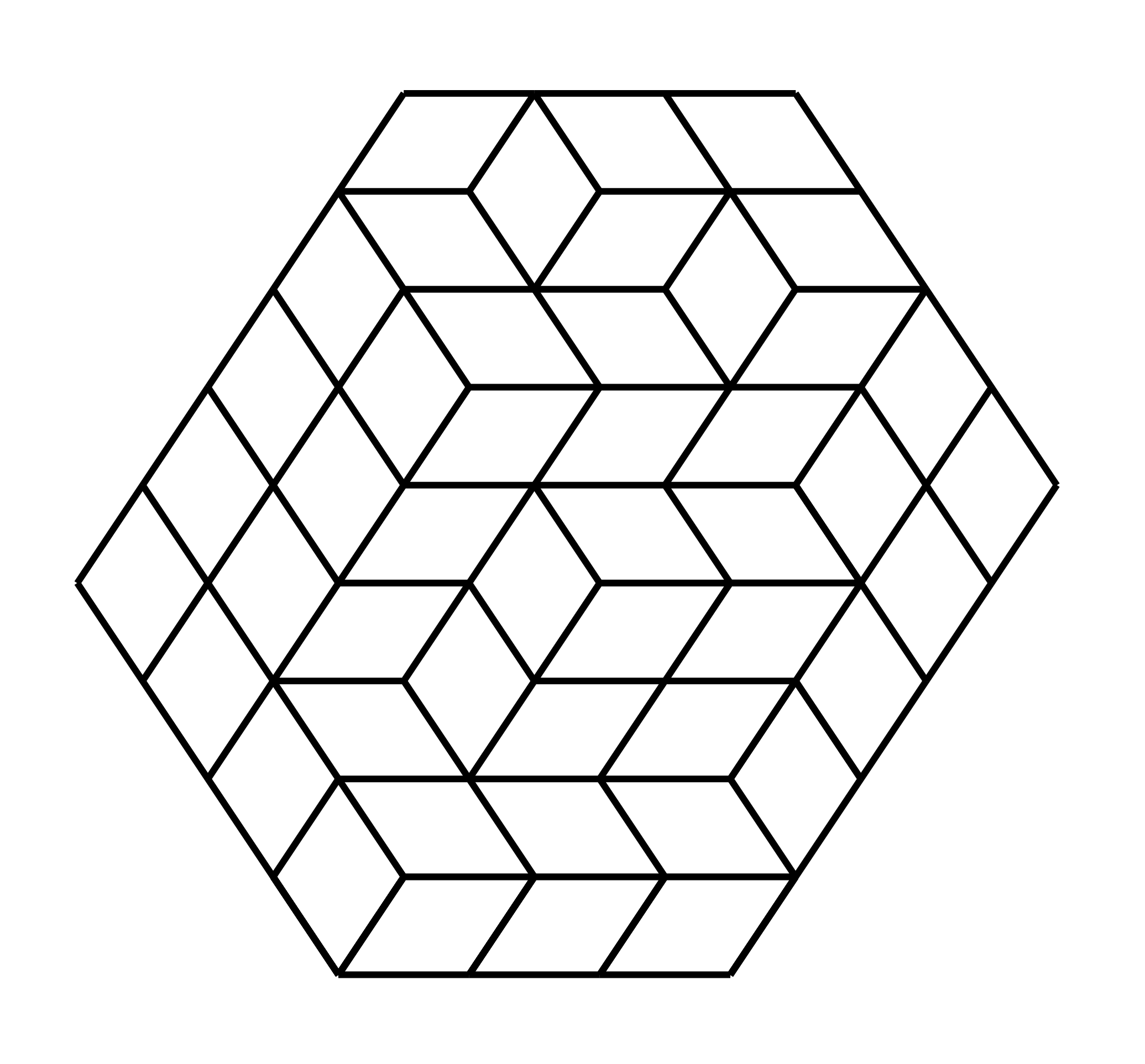}}
\hspace{1mm}
\subfigure[]
{\label{Mpath}\includegraphics[height=0.29\textwidth]{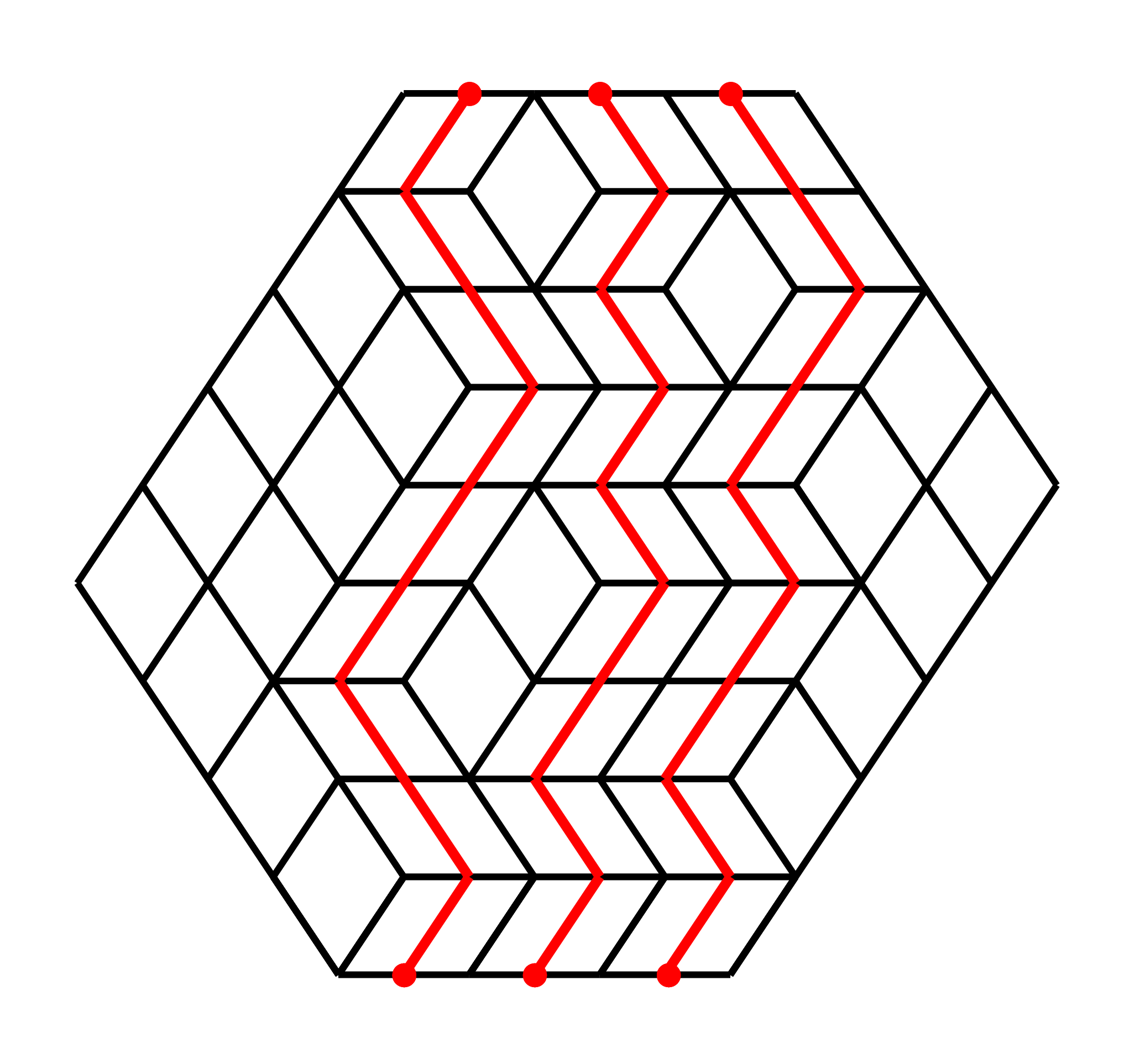}}
\caption{(a) An example of the plane partition represented as a pile of unit cubes in the $3 \times 4 \times 5$ box. (b) The corresponding lozenge tilings of the hexagon with side lengths $3,4,5,3,4,5$. (c) The corresponding family of non-intersecting lattice paths.}
\label{MacMahon}
\end{figure}

One can view (due to~\cite{DT89}) the plane partitions in an $a \times b \times c$ box as the lozenge tilings of a hexagonal region with side lengths $a,b,c,a,b,c$ (in cyclic order), see Figure~\ref{Mbox} and Figure~\ref{Mtile}. The number of these tilings is given by MacMahon's box formula~\cite[Sections 429 and 494]{Mac}
\begin{equation}
\prod_{i=1}^{a}\prod_{j=1}^{b}\prod_{k=1}^{c} \frac{i+j+k-1}{i+j+k-2}.
\end{equation}
Using the path method and that the set of tilings are in bijection with a family of non-intersecting paths on a square lattice (see Figure~\ref{Mpath}), the formula follows by evaluating the determinant of a matrix whose entries are binomial coefficients (see~\cite[Theorem 15]{GV89},\cite[Theorem 6.6]{Kra90}).

\begin{figure}[hbt!]
\centering
\subfigure[]
{\label{hexeven1}\includegraphics[width=0.47\textwidth]{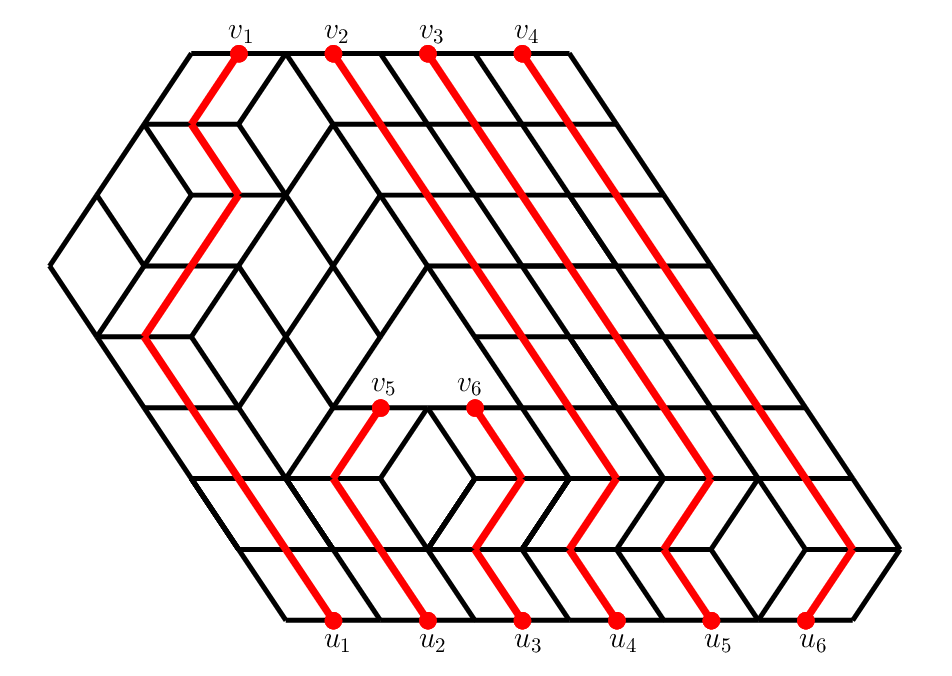}}
\hspace{1mm}
\subfigure[]
{\label{hexeven2}\includegraphics[width=0.47\textwidth]{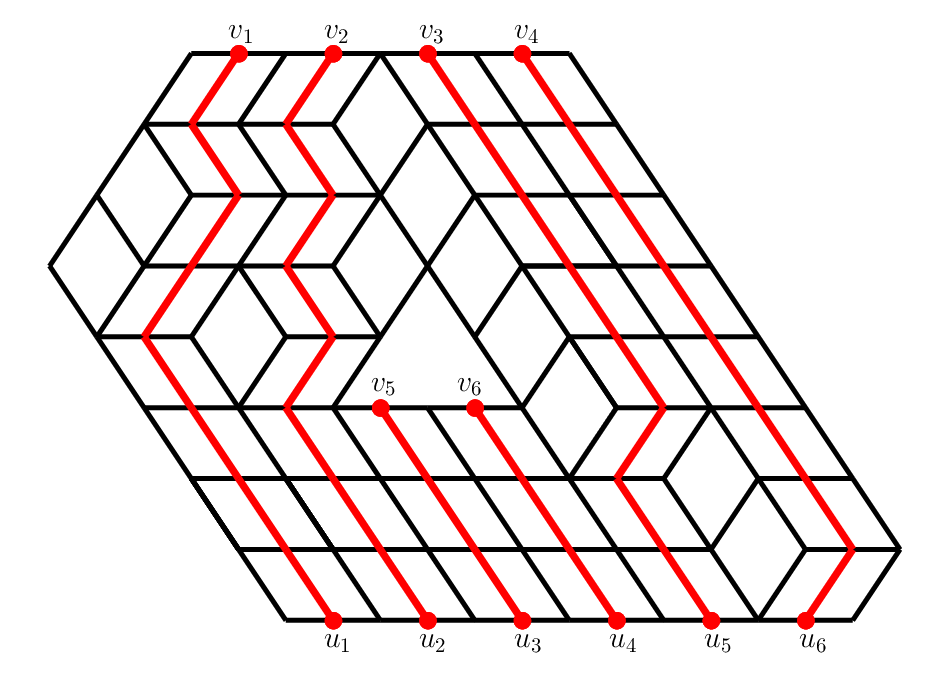}}
\caption{An example of two lozenge tilings of the hexagon with the even triangular hole and the corresponding lattice paths. The permutation (connection type) induced from these paths is $(1)(25364)$ in Figure~\ref{hexeven1} and is $(1)(2)(35)(46)$ in Figure~\ref{hexeven2}, they have the same permutation sign.}
\label{hexevenhole}
\end{figure}

\begin{figure}[hbt!]
\centering
\subfigure[]
{\label{hexodd1}\includegraphics[width=0.47\textwidth]{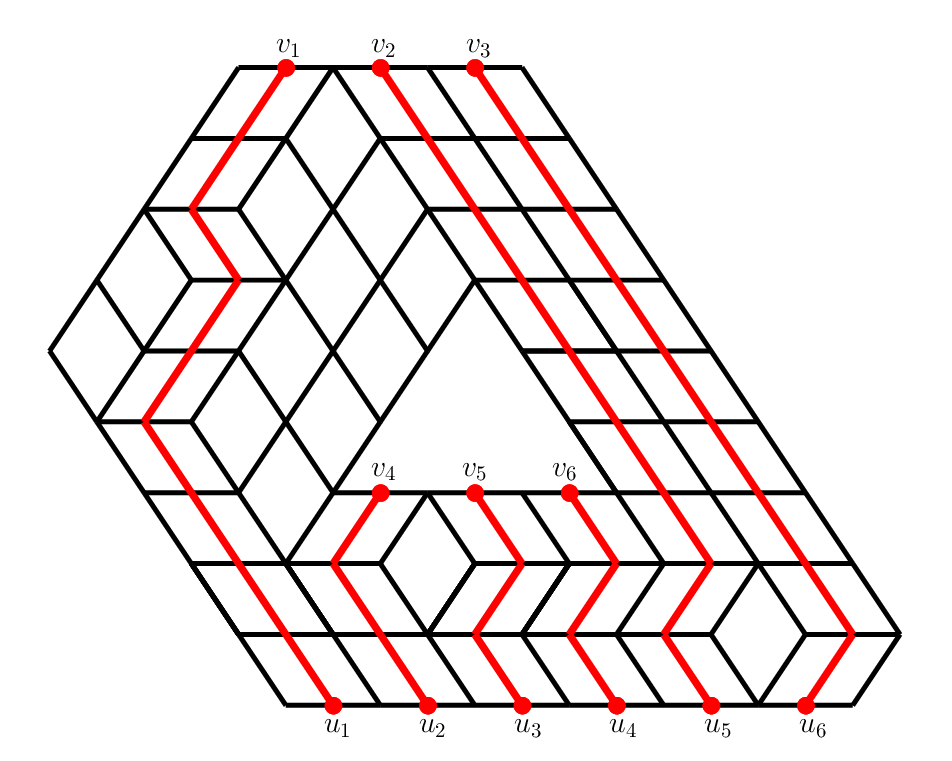}}
\hspace{1mm}
\subfigure[]
{\label{hexodd2}\includegraphics[width=0.47\textwidth]{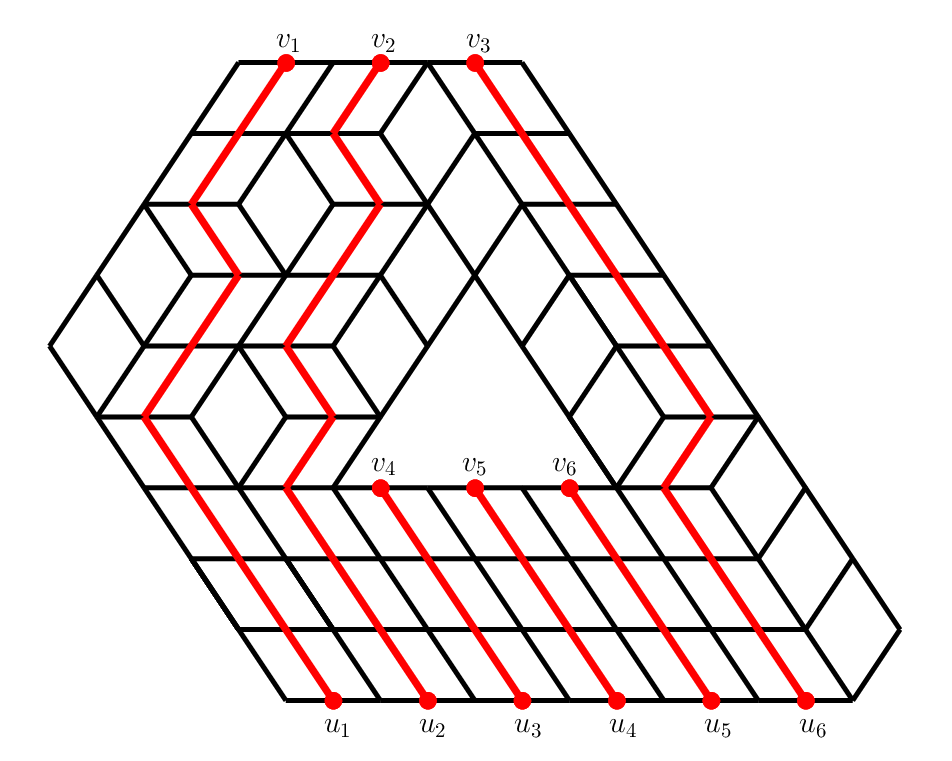}}
\caption{An example of two lozenge tilings of the hexagon with the odd triangular hole and the corresponding lattice paths. The permutation (connection type) induced from these paths is $(1)(24635)$ in Figure~\ref{hexodd1} and is $(1)(2)(3456)$ in Figure~\ref{hexodd2}, they have the opposite permutation sign.}
\label{hexoddhole}
\end{figure}

The path method~\cite[Section 3.1]{Propp15} can deal with the hexagonal region with triangular holes by adding suitable starting and ending points for the paths. The set of starting points and the set of ending points are not compatible in general (see Figure~\ref{hexevenhole} and Figure~\ref{hexoddhole}), so the left hand side of \eqref{eq.GV} does not reduce to a single term. This formula gives then a signed enumeration, not the straight enumeration of tilings. The signed enumeration of lozenge tilings was discussed in the paper by Cook II and Nagel \cite{CN17}, they defined the signed lozenge tilings in two natural ways. The first one comes from the lattice paths viewpoint, it is given by the determinant in the Lindstr{\"o}m-Gessel-Viennot theorem (as mentioned above). The other one comes from the perfect matchings viewpoint, it is given by the determinant of the bi-adjacency matrix of the corresponding bipartite graph. They showed that these two signed enumerations are equivalent. The signed lozenge tilngs was also mentioned in~\cite[Section 3]{Gil17}. The main result of this paper is to provide a determinant that gives the straight enumeration.

For the tilings on a hexagonal region with even triangular holes (triangular holes with even side length), it turns out that all the permutations induced from the family of lattice paths on the left hand side of~\eqref{eq.GV} have the same sign (see Figure~\ref{hexevenhole}), so Theorem~\ref{thm.GV1} does give the correct enumeration (see~\cite[Section 5]{CEKZ01},\cite[Section 3]{CK17},\cite{Gil18}, \cite{Con20} and \cite{Lai20}).

However, a hexagonal region containing an odd triangular hole (a triangular hole with odd side length) does not have this nice property (see Figure~\ref{hexoddhole}). In spite of this, it turns out that one can associate signs to each path so that the modified determinant gives the correct count. This idea is due to Krattenthaler~\cite{CKPrivate}, who showed how to deal with the odd triangular holes in a hexagonal region. This corresponds to families of non-intersecting lattice paths on the directed grid graph with compatible starting and ending points, except for a block of contiguous ending points in the middle. Our result extends his idea to acyclic digraphs with arbitrary starting and ending points.


The rest of this paper is organized as follows. In Section~$2$, we describe the specific way we associate signs to paths and state the main theorem (Theorem~\ref{thm.main}). In Section~$3$, we reduce the proof of the main theorem to two lemmas (Lemma~\ref{keylemma1} and Lemma~\ref{keylemma2}) and provide the proof of Lemma~\ref{keylemma1}. In Section~$4$, we define the transversal intersection number of paths and present the proof of Lemma~\ref{keylemma2}. We give some enumerative results in Section~$5$. The domino tilings of the mixed Aztec rectangle with arbitrary holes are discussed in Section~$6$.

\section{Statement of main results}

Drawing a graph in the plane in a meaningful way can be hard in general. The field of graph drawing deals with this problem. Algorithmic processes for visualizing graphs have been investigated by both mathematicians and computer scientists. We refer the interested reader to~\cite{BGTRI98} for more details.

An \textit{upward planar drawing} of a digraph $G$ is a drawing of $G$ on the Euclidean plane such that
\begin{itemize}
  \item each edge is drawn as a line segment that is either horizontal or up-pointing, and
  \item no two edges may intersect except at vertices of $G$.
\end{itemize}
In other words, if $e$ is a directed edge from $u$ to $v$ in an upward planar drawing graph, then the $y$-coordinate of $v$ is greater than or equal to the $y$-coordinate of $u$.

\begin{figure}[hbt!]
\centering
\subfigure[]
{\label{Fig1a}\includegraphics[height=0.23\textwidth]{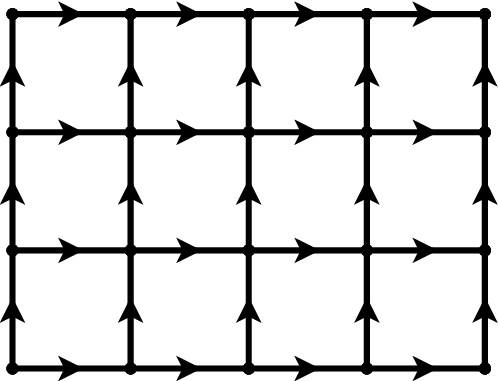}}
\hspace{6mm}
\subfigure[]
{\label{Fig1b}\includegraphics[height=0.23\textwidth]{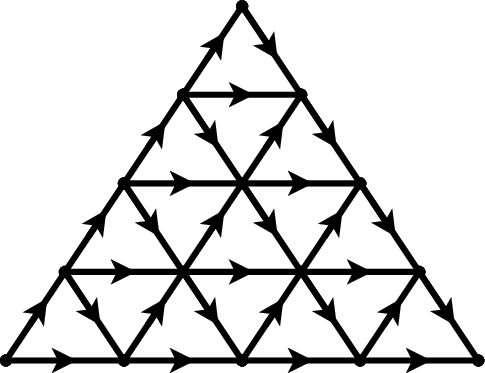}}
\hspace{6mm}
\subfigure[]
{\label{Fig1c}\includegraphics[height=0.23\textwidth]{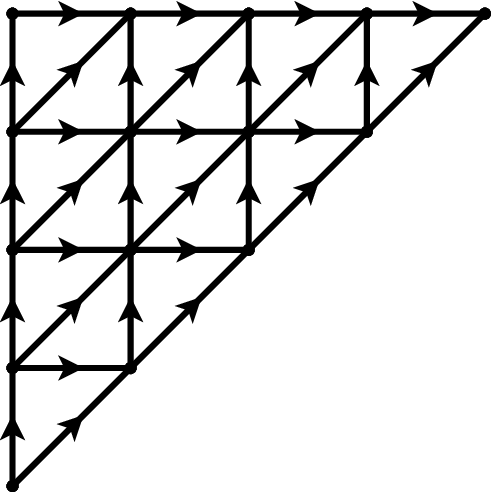}}
\caption{(a) The square lattice and the orientations. (b) The triangular lattice and the orientations. (c) An upward planar drawing of the triangular lattice in Figure~\ref{Fig1b}.}
\label{Fig1}
\end{figure}

Note that a square lattice (see Figure~\ref{Fig1a}) with horizontal (resp., vertical) edges oriented east (resp., north) is an upward planar drawing. A triangular lattice with orientations described in Figure~\ref{Fig1b} has the upward planar drawing shown in Figure~\ref{Fig1c}.

An $st\textit{-planar graph}$ is a planar, acyclic digraph with one source and one sink\footnote{A source is a vertex with no incoming edges. A sink is a vertex with no outgoing edges.}, so that these two special vertices lie on the outer face of the graph. The graphs that have an upward planar drawing are characterized by the following theorem stated in~\cite[Theorem 6.1]{BGTRI98}.
\begin{theorem}[Di Battista et al.~\cite{BGTRI98}]\label{thm.up+st}
  A graph $G$ has an upward planar drawing if and only if $G$ is a subgraph of an $st$-planar graph on the same vertex set.
\end{theorem}
Given an $st$-planar graph $\widetilde{G}$, consider the subgraph $G$ and directed paths on $G$ with given starting and ending points; these will be called \textit{marked points}. Let $U = \{u_1,\dotsc,u_n\}$ be the set of $n$ distinct starting points and $V = \{v_1,\dotsc,v_n\}$ be the set of $n$ distinct ending points. Now, we formally introduce our notations below.
\begin{itemize}
  \item $\mathscr{P}(u_i,v_j)$ denotes the set of paths going from $u_i \in U$ to $v_j \in V$.
  \item $\mathscr{P}^{\pi}(U,V)$ denotes the set of $n$-tuples of paths $(p_1,\dotsc,p_n)$, where $p_i \in \mathscr{P}(u_i,v_{\pi(i)})$ for $1 \leq i \leq n$. The permutation $\pi$ is called the \textit{connection type}.
  \item $\mathscr{P}(U,V)$ is the set of all $n$-tuples of paths connecting $U$ to $V$. In other words, $\mathscr{P}(U,V)$ is the union of $\mathscr{P}^{\pi}(U,V)$ over all the permutations $\pi \in \mathfrak{S}_n$.
  \item $\mathscr{P}_0(U,V)$ (resp., $\mathscr{P}^{\pi}_0(U,V)$) is the subset of $\mathscr{P}(U,V)$ (resp., $\mathscr{P}^{\pi}(U,V)$) consisting of non-intersecting $n$-tuples of paths.
\end{itemize}
\begin{definition}\label{def.extension}
  Let $s$ be the source and $t$ be the sink of the $st$-planar graph $\widetilde{G}$, and let $p \in \mathscr{P}(u,v)$ be a path in the subgraph $G$. The \textit{left side of the path $p$} is the closed region of the plane bounded by the following paths in $\widetilde{G}$:
  \begin{itemize}
    \item the leftmost path\footnote{The path obtained by taking the leftmost step at each stage.} from $s$ to $u$,
    \item the path $p$ itself,
    \item the leftmost path from $v$ to $t$, and
    \item the left boundary of $\widetilde{G}$ going from $s$ to $t$.
  \end{itemize}
  We write $L(p)$ for the collection of marked points of $U \cup V$ which are on the left side of the path $p$; this includes the starting point $u$ and the ending point $v$ of the path $p$.
\end{definition}
\begin{definition}\label{def.sign}
  The \textit{path sign} of a path $p \in \mathscr{P}(u,v)$ is defined to be
  \begin{equation}\label{eq.sign1}
    \sgn(p) = (-1)^{\left| L(p) \right|}.
  \end{equation}
  The path sign of an $n$-tuple of paths $P = (p_1,\dotsc,p_n) \in \mathscr{P}(U,V)$ is defined to be the product of all path signs of the $p_i$'s:
  \begin{equation}\label{eq.sign2}
    \sgn(P)=\prod_{i=1}^{n}\sgn(p_i).
  \end{equation}
\end{definition}
\begin{figure}[hbt!]
\centering
\subfigure[]
{\label{fig.ext1}\includegraphics[width=0.45\textwidth]{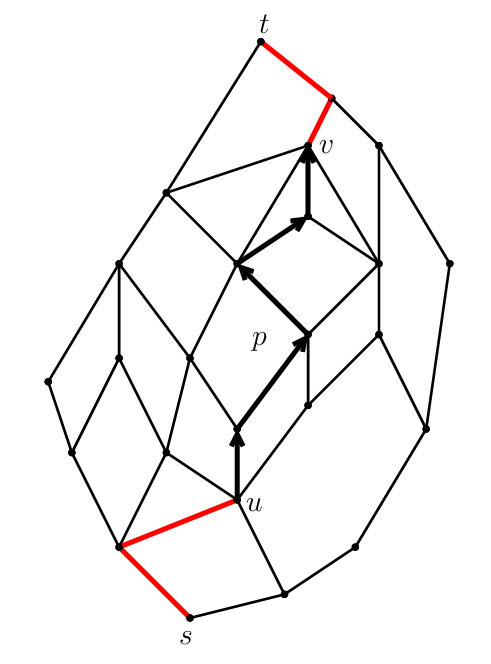}}
\hspace{10mm}
\subfigure[]
{\label{fig.ext2}\includegraphics[width=0.45\textwidth]{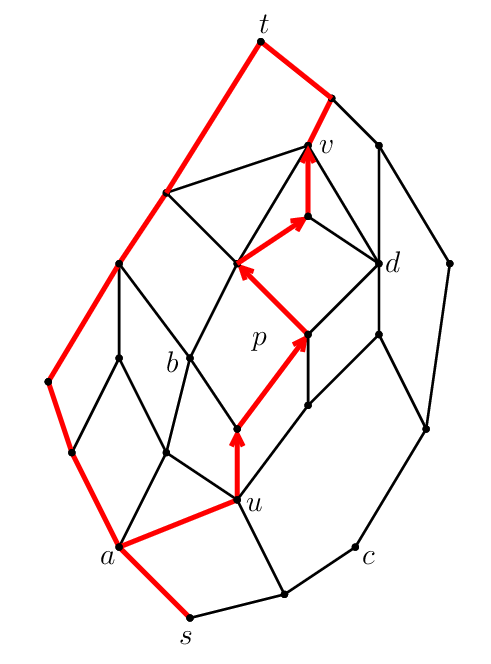}}
\caption{(a) The leftmost paths (red edges) from $s$ to $u$ and from $v$ to $t$. (b) The left side of the path $p$ is the region bounded by red edges.}
\label{fig.ext}
\end{figure}

Figure~\ref{fig.ext} shows an example of a graph $G$ having an upward planar drawing: the source and the sink are denoted by $s$ and $t$, respectively. In Figure~\ref{fig.ext1}, the leftmost paths from $s$ to $u$ and from $v$ to $t$ are drawn in red edges. The left side of $p$ is the region enclosed by the red edges. The four other marked points $a,b,c,d \in U \cup V$ (besides $u$ and $v$) are shown in Figure~\ref{fig.ext2}. We have $c,d \notin L(p)$ and $L(p)=\{a,b,u,v\}$. According to~\eqref{eq.sign1}, $\sgn(p) = (-1)^{4} = 1$.

Now, we are ready to state the main theorem.
\begin{theorem}\label{thm.main}
 Given an $st$-planar graph $\widetilde{G}$ and a subgraph $G$. Let $U=\{u_1,u_2,\dots,u_n\}$ and $V=\{v_1,v_2,\dots,v_n\}$ be two sets of $n$ distinct marked points of $G$. Let $M$ be the $n \times n$ matrix whose $(i,j)$-entry is
 \begin{equation}\label{M}
   \sum_{p \in \mathscr{P}(u_i,v_j)} \sgn(p)\wt(p).
 \end{equation}
 Then the total weight of families of non-intersecting paths connecting $U$ to $V$ is given by
 \begin{equation}\label{eq.main}
   GF(\mathscr{P}_{0}(U,V)) = \sum_{\pi \in \mathfrak{S}_n} GF(\mathscr{P}^{\pi}_{0}(U,V)) = \left| \det M \right|.
 \end{equation}
\end{theorem}
\begin{remark}
  On a directed acyclic graph, the Lindstr{\"o}m-Gessel-Viennot theorem gives the signed enumeration \eqref{eq.GV} or the simplified case \eqref{eq.cor} when the starting and ending points are compatible. Our result (Theorem \ref{thm.main}) gives the straight enumeration for arbitrary (equinumerous) starting and ending points for any digraph that has an upward planar drawing.
\end{remark}
\begin{example}
Figure~\ref{fig.example} shows the $6 \times 6$ grid graph, with the horizontal edges oriented east, weighted by $x$, and vertical edges oriented north, weighted by $y$. In this example $n=2$, so we have two starting points $U=\{u_1,u_2\}$ and two ending points $V=\{v_1,v_2\}$. We would like to find the weighted sum of non-intersecting paths from $U$ to $V$.

Note that there are two connection types: $\mathsf{id}$ and $(12)$ (see Figure~\ref{fig.ex1} and Figure~\ref{fig.ex2}, respectively). Thus the Lindstr{\"o}m-Gessel-Viennot theorem gives the signed weighted count
$$\sgn(\mathsf{id})GF\left(\mathscr{P}_{0}^{\mathsf{id}}(U,V)\right) + \sgn((12))GF\left( \mathscr{P}_{0}^{(12)}(U,V)\right).$$
By contrast, Theorem~\ref{thm.main}, gives the straight weighted count of families of non-intersecting paths:
$$GF\left(\mathscr{P}_{0}^{\mathsf{id}}(U,V)\right) + GF\left(\mathscr{P}_{0}^{(12)}(U,V)\right).$$
\begin{figure}[hbt!]
    \centering
    \subfigure[The paths with connection type $\mathsf{id}$, $p,q$ are two different types of paths from $u_2$ to $v_2$.]
    {\label{fig.ex1}\includegraphics[width=0.4\textwidth]{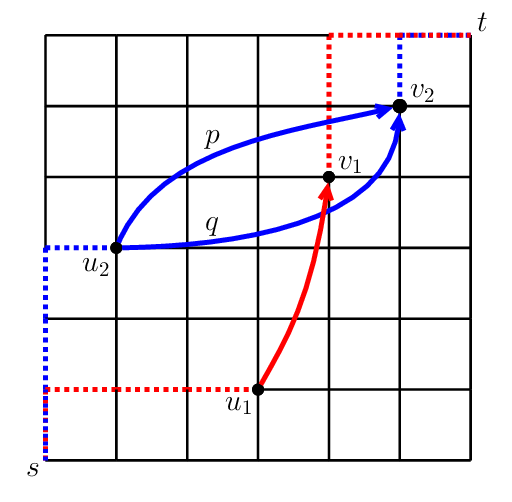}}
    \hspace{10mm}
    \subfigure[The paths with connection type $(12)$, $p^{\prime},q^{\prime}$ are two different types of paths from $u_1$ to $v_2$.]
    {\label{fig.ex2}\includegraphics[width=0.4\textwidth]{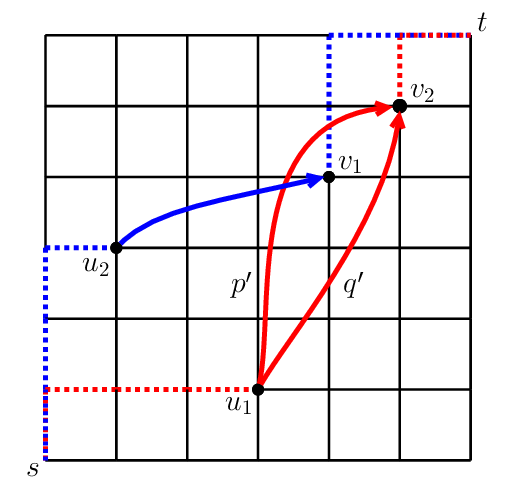}}
    \caption{The $6 \times 6$ square lattice with two starting and two ending points.}
    \label{fig.example}
  \end{figure}

In order to get the explicit expression, note that in Figure~\ref{fig.example} the dotted edges form the leftmost paths connecting $s$ to $u_i$ and $v_{\pi(i)}$ to $t$. In Figure~\ref{fig.ex1}, there are $\binom{6}{4}$ paths from $u_2$ to $v_2$. Of these, $\binom{4}{2}$ of them (exemplified by $p$) have only $u_2$ and $v_2$ on their left sides, while $9$ of them (exemplified by $q$, including the paths passing through $v_1$) have $u_2$, $v_2$ and $v_1$ on their left sides. After considering the weights and the path signs, the $(2,2)$-entry of the matrix $M$ works out to be $6x^4y^2-9x^4y^2 = -3x^4y^2$. On the other hand, there are $\binom{4}{1}$ paths joining $u_1$ and $v_1$, and for all of them the left sides contain $u_1$, $v_1$ and $u_2$; thus, they all have path sign $-1$. The $(1,1)$-entry of the matrix $M$ is therefore $-4xy^3$.

Similarly, in Figure~\ref{fig.ex2}, there are $\binom{4}{3}$ paths joining $u_2$ and $v_1$, and $u_2$ and $v_1$ lie on their left sides. Thus the path signs of these are positive, and the $(2,1)$-entry of the matrix $M$ is $4xy^3$. Finally, there are $\binom{6}{2}$ paths from $u_1$ to $v_2$: $14$ of them (illustrated by $q^{\prime}$) have $v_1$ contained in their left sides, so all the four points are on their left sides; this gives positive path signs. Only $1$ of them (represented by $p^{\prime}$) does not have $v_1$ on its left side; its path sign is negative. The $(1,2)$-entry of the matrix $M$ is therefore $14x^2y^4-x^2y^4 = 13x^2y^4$.

Combining the above discussions and Theorem~\ref{thm.main}, we obtain that the straight enumeration of non-intersecting paths connecting $U$ to $V$ is equal to
$$ \left| \det M \right| =
\left| \det \begin{pmatrix}
-4xy^3 & 13x^2y^4 \\
4x^3y &  -3x^4y^2
\end{pmatrix}\right|
= 40x^5y^5.
$$
\end{example}

\section{Two lemmas}

In this section we present two lemmas (Lemma~\ref{keylemma1} and Lemma~\ref{keylemma2}), and show how the proof of Theorem~\ref{thm.main} can be readily obtained as a consequence. We prove the former in this section, but postpone the proof of the latter to Section $4$.

The proof of the Lindstr{\"o}m-Gessel-Viennot theorem contains a clever sign-reversing and weight-preserving involution $\phi$ on families of paths with given starting points $U$ and ending points $V$; see for instance~\cite[Theorem 1.2]{Stem90} and~\cite[Section 5.4]{Aig07}. This involution is described below.

First, choose a total order of the vertices of the graph. Consider any $n$-tuple of intersecting paths with the given connection type $\pi$, say $P=(p_1,\dotsc,p_n) \in \mathscr{P}^{\pi}(U,V)$. Find the least vertex $w$ which is the vertex of intersection of paths, if there are more than two paths meet at $w$, then select two distinct paths $p_i$ and $p_j$, $i \neq j$, with the smallest indices. This makes $w$, $p_i$ and $p_j$ unique for each $n$-tuple of intersecting paths.

Second, create the new path $p_i^{*}$ (resp., $p_j^{*}$) by concatenating the first half of $p_i$ (resp., $p_j$) up to $w$ and the second half of $p_j$ (resp., $p_i$) after the vertex $w$. Note that the ending points of these two paths interchanged. Replace $p_i$ by $p_i^{*}$ and $p_j$ by $p_j^{*}$, then we obtain the new family of intersecting paths
\begin{equation}\label{eq.involution}
  P^{*}=(p_1,\dots,p_{i-1},p_i^{*},p_{i+1},\dots,p_{j-1},p_j^{*},p_{j+1},\dots,p_n),
\end{equation}
which has the connection type $\sigma \pi$, where $\sigma = (\pi(i)\pi(j))$.

On the other hand, the $n$-tuple of non-intersecting paths is defined to be the fixed point of this involution. So, the involution $\phi: \mathscr{P}(U,V) \rightarrow \mathscr{P}(U,V)$ is given by
\begin{equation}\label{eq.involution0}
  \phi(P) = \left\{
  \begin{array}{ll}
  P^{*}, & \text{if $P$ is a family of intersecting paths.} \\
  P, & \text{if $P$ is a family of non-intersecting paths.}
  \end{array}
  \right.
\end{equation}

It is well-known that this involution is sign-reversing (for permutations) and weight-preserving. Our first lemma, below, show that this involution also preserves the path sign.
\begin{lemma}\label{keylemma1}
  Let $P=(p_1,\dots,p_n)$ be a family of paths. Then the image $\phi(P)$ of $P$ under the involution \eqref{eq.involution0} satisfies
  \begin{equation}\label{eq.sign+involution}
    \sgn(P) = \sgn(\phi(P)).
  \end{equation}
\end{lemma}
\begin{proof}
  It is trivial that \eqref{eq.sign+involution} holds for families of non-intersecting paths. Now, we assume $P$ is a family of intersecting paths.

  Figure~\ref{fig.keylemma1-1} illustrates two paths $p_i$ and $p_j$ which meet at the vertex $x$, their corresponding images under the involution are given by $p_i^{*}$ and $p_j^{*}$ in Figure~\ref{fig.keylemma1-2}. The left side of $p_i$ and $p_i^{*}$ (resp., $p_j$ and $p_j^{*}$) is the region bounded by red (resp., blue) curves. It is easy to see that the union of the left side of $p_i$ and the left side of $p_j$ equals the union of the left side of $p_i^{*}$ and the left side of $p_j^{*}$. As a consequence,
  \begin{equation}\label{eq.union}
    \left|L(p_i)\right| + \left|L(p_j)\right| = \left|L(p_i^{*})\right| + \left|L(p_j^{*})\right|.
  \end{equation}

  By \eqref{eq.union}, we obtain the identity of the path signs
  \begin{equation*}
    \sgn(p_i)\sgn(p_j) = (-1)^{\left| L(p_i) \right| + \left| L(p_j) \right|} = (-1)^{\left| L(p_i^{*}) \right| + \left| L(p_j^{*}) \right|} = \sgn(p_i^{*})\sgn(p_j^{*}).
  \end{equation*}
  Therefore,
  \begin{equation*}
    \sgn(P) = \left( \prod_{k \neq i,j} \sgn(p_k) \right) \sgn(p_i)\sgn(p_j) = \left( \prod_{k \neq i,j} \sgn(p_k) \right) \sgn(p_i^{*})\sgn(p_j^{*}) = \sgn(P^{*}).
  \end{equation*}
  \begin{figure}[hbt!]
    \centering
    \subfigure[]
    {\label{fig.keylemma1-1}\includegraphics[width=0.42\textwidth]{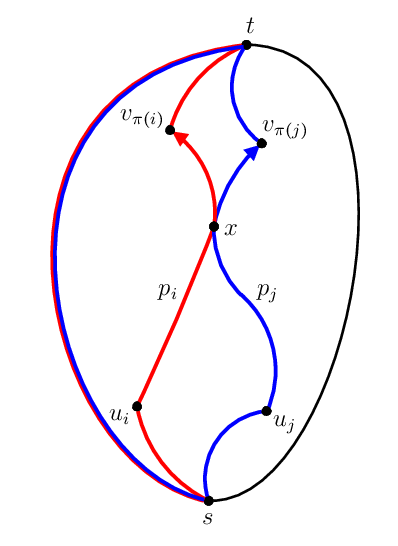}}
    \hspace{5mm}
    \subfigure[]
    {\label{fig.keylemma1-2}\includegraphics[width=0.42\textwidth]{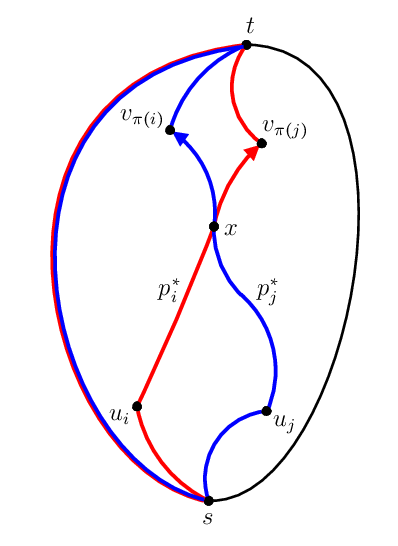}}
    \caption{(a) Illustration of two paths $p_i$ (red) and $p_j$ (blue), and their left sides. (b) The image of $p_i$ and $p_j$ under the involution, and their left sides.}
    \label{fig.keylemma1}
  \end{figure}
\end{proof}
\begin{lemma}\label{keylemma2}
  Given a connection type $\pi$, if the families of non-intersecting paths $P \in \mathscr{P}^{\pi}_{0}(U,V)$ and $Q \in \mathscr{P}^{\sigma \pi}_{0}(U,V)$ exist for some transposition $\sigma \in \mathfrak{S}_n$, then
  \begin{equation}\label{sign+non-intersetion1}
    \sgn(P) = -\sgn(Q).
  \end{equation}
  Therefore,
  \begin{equation}\label{eq3}
    \sgn(\pi)\sgn(P) = -\sgn(\pi)\sgn(Q) = \sgn(\sigma \pi)\sgn(Q).
  \end{equation}
In other words, as $P$ runs over $\mathscr{P}_0(U,V)$, if $\pi_P$ is the connection type of $P$, then $\sgn(\pi_P)\sgn(P)$ is constant.
\end{lemma}

The proof of Lemma \ref{keylemma2} requires some new ideas and will be presented in the next section. Now, we are able to prove our main theorem.
\begin{proof}[Proof of Theorem~\ref{thm.main}]
We remind the reader that the $(i,j)$-entry of the matrix $M$ is given by $\sum_{p \in \mathscr{P}(u_i,v_j)} \sgn(p)\wt(p)$. The determinant of $M$ can be interpreted as the following due to the independence of the indices $\pi,p_1,\dotsc,p_n$:
\begin{align}
  \det M & = \sum_{\pi \in \mathfrak{S}_n} \sgn(\pi) \left( \sum_{p_1 \in \mathscr{P} \left( u_1,v_{\pi(1)} \right) }\sgn(p_1)\wt(p_1) \right) \cdots \left( \sum_{p_n \in \mathscr{P} \left( u_n,v_{\pi(n)} \right) }\sgn(p_n)\wt(p_n) \right) \nonumber \\
  & = \sum_{\pi \in \mathfrak{S}_n} \sgn(\pi) \left( \sum_{P=(p_1,\dotsc,p_n) \in \mathscr{P}^{\pi}(U,V)} \sgn(P)\wt(P) \right) \nonumber \\
  & = \sum_{(\pi, P) \in \mathfrak{S}_n \times \mathscr{P}(U,V)} \sgn(\pi)\sgn(P)\wt(P), \label{eq1}
\end{align}
where the last summation is over all pairs $(\pi,P)$ so that $P$ has connection type $\pi$.

The involution $\phi$ described above induces $\overline{\phi}: \mathfrak{S}_n \times \mathscr{P}(U,V) \rightarrow \mathfrak{S}_n \times \mathscr{P}(U,V)$,
\begin{equation}\label{eq.induceinv}
\overline{\phi}\left((\pi,P)\right)=
\left\{
\begin{array}{ll}
 (\sigma\pi,P^{*}), & \text{if $P$ is a family of intersecting paths} \\
 (\pi,P), & \text{if $P$ is a family of non-intersecting paths}
\end{array}
\right.
\end{equation}
(recall that $\sigma\pi$ is the connection type of $P^*$, and $\sigma$ is a transposition).
We apply $\overline{\phi}$ on $\mathfrak{S}_n \times \mathscr{P}(U,V)$. If $P$ is a family of intersecting paths, then $\overline{\phi}$ preserves $\wt(P)$ (by definition) and $\sgn(P)$ (by Lemma~\ref{keylemma1}) but reverses the permutation signs. As a consequence, all the families of intersecting paths are cancelled out in~\eqref{eq1}, the only contribution comes from families of non-intersecting paths in $\mathscr{P}_{0}(U,V)$. Thus~\eqref{eq1} can be simplified as
\begin{equation}
  \sum_{(\pi, P_0) \in \mathfrak{S}_n \times \mathscr{P}_{0}(U,V)} \sgn(\pi)\sgn(P_0)\wt(P_0), \label{eq2}
\end{equation}
where the summation is over all pairs $(\pi,P_0)$ so that $P_0$ has connection type $\pi$.

Using Lemma~\ref{keylemma2}, \eqref{eq2} can be written as
\begin{equation}\label{eq4}
  \pm \sum_{\pi \in \mathfrak{S}_n} \left(\sum_{P_0 \in \mathscr{P}^{\pi}_{0}(U,V)} \wt(P_0) \right) = \pm \sum_{\pi \in \mathfrak{S}_n} GF(\mathscr{P}^{\pi}_{0}(U,V)) = \pm GF(\mathscr{P}_0(U,V)).
\end{equation}
This completes the proof of our main theorem.
\end{proof}
\begin{remark}
  The reader might notice that an upward planar drawing of a graph $G$ may not be unique. If we work on a graph $G$ with two different drawings, then the number of marked points in $U \cup V$ on the left side of a given path $p$ may be different in the two drawings; this may yield different path signs for $p$, and hence different matrices $M$ in Theorem \ref{thm.main}. However, both drawings must give the same result, so their determinants must be equal. We could obtain this way either some new determinant identities, or new proofs of some known determinant identities by constructing suitable drawings of the same graph. We leave this direction to be pursued by the interested reader.
\end{remark}

\section{Intersection number of a family of paths}

In this section we give the proof of Lemma \ref{keylemma2}. Since Lemma \ref{keylemma2} involves only families of non-intersecting paths, we assume throughout this section that the paths $p \in \mathscr{P}(u_i,v_j)$ we work with do not contain any marked points of $U \cup V$ other than the starting and ending points of $p$.

Consider a family of paths $P = (p_1,\dotsc,p_n)\in \mathscr{P}(U,V)$ connecting $U$ to $V$. For $i \neq j$, let $p_i \cap p_j = \{c_1, \dotsc, c_{\ell} \}$ be the (possibly empty) collection of maximal subpaths $c_k$ that $p_i$ and $p_j$ share; this includes subpaths of length zero, that is, single vertices. It is clear that the paths $p_i$ and $p_j$ are non-intersecting if and only if $p_i \cap p_j = \emptyset$.

\begin{definition}\label{def.intersection}
  We say that $c_{k}$ is a \textit{transversal intersection} of $p_i$ and $p_j$, $i \neq j$, if
  \begin{itemize}
    \item $c_{k} \in p_i \cap p_j$, and
    \item the path $p_i$ arrives and leaves $c_k$ from different sides of the path $p_j$.\footnote{It is easy to see that this is equivalent to saying that the path $p_j$ reaches and leaves $c_k$ from different sides of the path $p_i$; see Figure \ref{fig.keylemma1-2}.}
  \end{itemize}
  The \textit{(transversal) intersection number of $p_i$ and $p_j$}, denoted by $I(p_i,p_j)$, is defined to be the total number of transversal intersections in $p_i \cap p_j$. For a family of paths $P=(p_1,\dots,p_n)$, the \textit{(transversal) intersection number of $P$}, denoted by $I(P)$, is given by
  \begin{equation}\label{eq.mod2int}
    I(P) = \sum_{1 \leq i < j \leq n} I(p_i,p_j).
  \end{equation}
  In particular, if $P$ is a family of non-intersecting paths, then $I(P)=0$.
\end{definition}

Our proof of Lemma \ref{keylemma2} is based on the following results.
\begin{lemma}\label{lem.intersection}
  Let $P=(p_1,\dotsc,p_{i-1},p_i,p_{i+1},\dotsc,p_n)$ and $P^{\prime}=(p_1,\dotsc,p_{i-1},p_i^{\prime},p_{i+1},\dotsc,p_n)$ be two families of paths with the same connection type $\pi$, where $p_i$ and $p_i^{\prime}$ are two paths from $u_i$ to $v_{\pi(i)}$ with $|L(p_i)| -  |L(p^{\prime}_i)| = d$. Then
  \begin{equation}\label{eq.int}
    I(P^{\prime}) \equiv I(P) + d \pmod{2}.
  \end{equation}
\end{lemma}
\begin{proof}
The only difference of $P$ and $P^{\prime}$ is the $i$th path, by \eqref{eq.mod2int}, we have
\begin{equation}\label{eq.pp}
  I(P)-I(P^{\prime}) = \sum_{\substack{1 \leq m \leq n \\ m \neq i}} I(p_i,p_m) - \sum_{\substack{1 \leq m \leq n \\ m \neq i}} I(p_i^{\prime},p_m).
\end{equation}
Suppose that the paths $p_i$ and $p_i^{\prime}$ meet at vertices $x_0=u_i$, $x_1,\dotsc,x_k = v_{\pi(i)}$. The region\footnote{This region reduces to a path in case $p_i$ and $p^{\prime}_i$ coincide between $x_{\ell-1}$ and $x_{\ell}$.} bounded by the subpaths of $p_i$ and $p_i^{\prime}$ connecting $x_{\ell-1}$ to $x_{\ell}$ is denoted by $D_{\ell}$, for $1 \leq \ell \leq k$. Figure~\ref{fig.lemma4-21} shows an example when $k=3$. Let $L^{\mathsf{int}}$ be the collection of marked points of $U \cup V$ which lie in the intersection of the left side of $p_i$ and the left side of $p_i^{\prime}$.

For each $\ell = 1,\dotsc,k$, we assume that $d_{\ell}$ of the marked points of $U \cup V$ are in the interior of $D_{\ell}$ (there are no marked points on the boundary of $D_{\ell}$ by our assumption in this section). Then we have
\begin{equation}\label{eq.extension}
  \left| L(p_i) \right| + \left| L(p_i^{\prime}) \right| = d_1 + \dotsc + d_k + 2|L^{\mathsf{int}}|.
\end{equation}
By \eqref{eq.extension}, we obtain
\begin{align}
\sum_{\ell=1}^{k} d_{\ell} & \equiv \left| L(p_i) \right| + \left| L(p_i^{\prime}) \right| - 2|L^{\mathsf{int}}| \pmod{2} \nonumber \\
 & \equiv \left| L(p_i) \right| - \left| L(p_i^{\prime}) \right| \pmod{2} \nonumber \\
 & \equiv d \pmod{2}. \label{eq.extension1}
\end{align}

By definition, each marked point in the interior of $D_{\ell}$ is either a starting point or an ending point of a path. We divide the paths that intersect the boundary of $D_\ell$ in the following three types (illustrated in Figure~\ref{fig.lemma4-22}):
\begin{enumerate}
  \item Paths with the starting and the ending points contained in $D_{\ell}$ (these are shown in red in Figure~\ref{fig.lemma4-22}). Suppose there are $d_{\ell,1}$ such paths (accounting for $2d_{\ell,1}$ points in $D_{\ell}$). Each path of this type intersects transversally the boundary of $D_{\ell}$ an even number times. We write $\sum_{j=1}^{d_{\ell,1}}(2n_{j,1})$ for the total number of these transversal intersections, where $n_{j,1} \in \mathbb{Z}_{\geq 0}$.
  \item Paths having exactly one of the starting and ending point in $D_{\ell}$ (shown in blue). Suppose there are $d_{\ell,2}$ such paths. Each path of this type intersects transversally the boundary of $D_{\ell}$ an odd number times. Thus the total number of these transversal intersections can be written as $\sum_{j=1}^{d_{\ell,2}}(2n_{j,2}+1)$, where $n_{j,2} \in \mathbb{Z}_{\geq 0}$.
  \item Paths with no starting and ending point in $D_{\ell}$ (shown in green). Each such path intersects the boundary of $D_{\ell}$ an even number times. Let $\sum_{j=1}^{d_{\ell,3}}(2n_{j,3})$, $n_{j,3} \in \mathbb{Z}_{\geq 0}$ be the total number of these transversal intersections.
\end{enumerate}
\begin{figure}[hbt!]
    \centering
    \subfigure[]
    {\label{fig.lemma4-21}\includegraphics[height=0.57\textwidth]{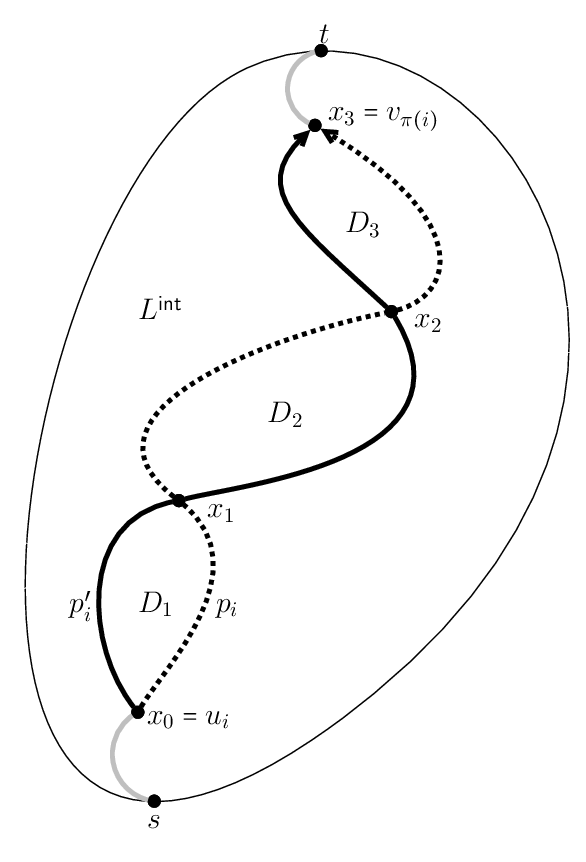}}
    \hspace{5mm}
    \subfigure[]
    {\label{fig.lemma4-22}\includegraphics[height=0.57\textwidth]{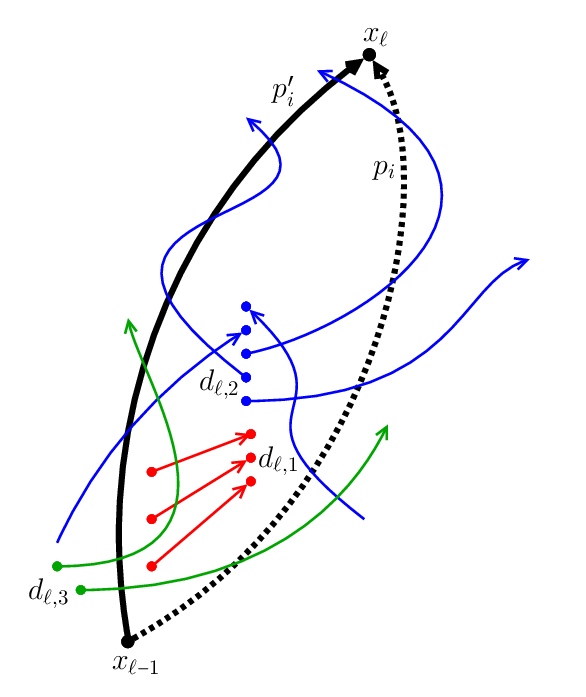}}
    \caption{(a) The paths $p_i$ (dotted curve) and $p_i^{\prime}$ (solid curve) meet at $x_0$, $x_1$, $x_2$ and $x_3$, they form three closed regions $D_1$, $D_2$ and $D_3$. $L^{\mathsf{int}}$ is the intersection of the left side of $p_i$ and the left side of $p_i^{\prime}$. (b) Illustration of the boundary of $D_{\ell}$ meets three different types (in different colors) of paths.}
    \label{fig.lemma4-2}
\end{figure}

For each $\ell = 1,\dotsc,k$, let $T_{\ell}$ be the total number of transversal intersections of the boundary of $D_\ell$ with paths $p_m$, for $1 \leq m \leq n, m \neq i$. Then we have
\begin{equation}\label{eq.intersections}
  T_{\ell} = \sum_{j=1}^{d_{\ell,1}}(2n_{j,1}) + \sum_{j=1}^{d_{\ell,2}}(2n_{j,2}+1) + \sum_{j=1}^{d_{\ell,3}}(2n_{j,3}),\ n_{j,1},n_{j,2},n_{j,3} \in \mathbb{Z}_{\geq 0}.
\end{equation}
By definition, we also have
\begin{equation}\label{eq.points}
  d_{\ell} = 2d_{\ell,1} + d_{\ell,2}.
\end{equation}
Combining \eqref{eq.intersections}, \eqref{eq.points} and summing over all $\ell$, we obtain
\begin{align}
\sum_{\ell=1}^{k} T_{\ell} & \equiv \sum_{\ell=1}^{k} d_{\ell,2} \pmod{2} \nonumber \\
 & \equiv \sum_{\ell=1}^{k} \left( d_{\ell} - 2d_{\ell,1} \right) \pmod{2} \nonumber \\
 & \equiv \sum_{\ell=1}^{k} d_{\ell} \pmod{2}. \label{eq.p+i}
\end{align}

We note that the union of the boundaries of $D_{\ell}$ consists of paths $p_i$ and $p_i^{\prime}$. So,
\begin{equation}\label{eq.totalint}
  \sum_{\ell=1}^{k} T_{\ell} = \sum_{\substack{1 \leq m \leq n \\ m \neq i}} I(p_i,p_m) + \sum_{\substack{1 \leq m \leq n \\ m \neq i}} I(p_i^{\prime},p_m)
\end{equation}
From the above discussion, we obtain
\begin{align*}
I(P) + I(P^{\prime}) & \equiv \sum_{\substack{1 \leq m \leq n \\ m \neq i}} I(p_i,p_m) + \sum_{\substack{1 \leq m \leq n \\ m \neq i}} I(p_i^{\prime},p_m) \pmod{2} & \text{(by \eqref{eq.pp})} \\
 & \equiv \sum_{\ell=1}^{k} T_{\ell} \pmod{2} & \text{(by \eqref{eq.totalint})}\\
 & \equiv \sum_{\ell=1}^{k} d_{\ell} \pmod{2} & \text{(by \eqref{eq.p+i})}\\
 & \equiv d \pmod{2}. & \text{(by \eqref{eq.extension1})}
\end{align*}
Therefore, $I(P^{\prime}) \equiv I(P) + d \pmod{2}$.
\end{proof}
\begin{lemma}\label{lem.sign+nonint}
Let $P \in \mathscr{P}^{\pi}(U,V)$ be a family of paths and $Q \in \mathscr{P}_{0}^{\pi}(U,V)$ be a family of non-intersecting paths having the same connection type $\pi$. Then
\begin{equation}\label{eq.sign+nonint}
  \sgn(P) = \sgn(Q)(-1)^{I(P)}.
\end{equation}
\end{lemma}
\begin{proof}
Let $P=(p_1,\dotsc,p_n)$, $Q = (q_1,\dotsc,q_n)$ and set $f_i=|L(p_i)| - |L(q_i)|$, for $i=1,2,\dotsc,n$. By Definition~\ref{def.sign}, we have $\sgn(p_i) = \sgn(q_i)(-1)^{f_i}$, for $1 \leq i \leq n$. Then
\begin{equation}\label{eq.01}
  \sgn(P) = \sgn(Q)(-1)^{f_1+\cdots+f_n}.
\end{equation}

We consider the sequence of families of paths $A_0 = Q, A_1,\dotsc,A_{n-1}, A_n = P$, where $A_i = (p_1,\dotsc,p_{i-1},p_i,q_{i+1},q_{i+2},\dotsc,q_n)$ for $1 \leq i \leq n-1$. Since $A_i$ and $A_{i-1}$ only differ in their $i$th paths, Lemma~\ref{lem.intersection} implies that
\begin{equation}\label{eq.02}
  I(A_i) \equiv I(A_{i-1}) + f_i \pmod{2},\ i=1,\dotsc,n.
\end{equation}
Therefore, we obtain
\begin{equation}\label{eq.03}
  I(P) \equiv I(Q) + f_1 + \cdots + f_n \pmod{2}.
\end{equation}

Since $I(Q)=0$, \eqref{eq.03} implies that $I(P)$ and $f_1 + \cdots + f_n$ have the same parity. Then the equation \eqref{eq.01} is equivalent to the equation \eqref{eq.sign+nonint}, as desired.
\end{proof}
%

We remind the reader that, in Lemma \ref{keylemma2}, if the connection types of two families of non-intersecting paths $P$ and $Q$ just differ from a transposition, then we would like to show that their path signs satisfy $\sgn(P) = -\sgn(Q)$.

\begin{proof}[Proof of Lemma \ref{keylemma2}]
Write the transposition $\sigma$ in the statement as $\sigma = (\pi(i)\pi(j))$, where $1 \leq i \neq j \leq n$. Let $P = (p_1,\dotsc,p_n) \in \mathscr{P}^{\pi}_0(U,V)$ and $Q = (q_1,\dotsc,q_n) \in \mathscr{P}^{\sigma\pi}_0(U,V)$ be the families of non-intersecting paths with the connection types $\pi$ and $\sigma\pi$, respectively.

Let $\mathbf{y}(v)$ be the $y$-coordinate of the vertex $v$ of $\widetilde{G}$. Since the families of paths $P$ and $Q$ exist, it follows in particular that there are paths in $G$ (recall that $G$ is a subgraph of $\widetilde{G}$ having the same vertex set) from $u_i$ to $v_{\pi(i)}$, from $u_i$ to $v_{\pi(j)}$, from $u_j$ to $v_{\pi(i)}$ and from $u_j$ to $v_{\pi(j)}$. Therefore, since $G$ is drawn in an upward-planar fashion, $\mathbf{y}(v_{\pi(i)}) \geq \mathbf{y}(u_i)$, $\mathbf{y}(v_{\pi(j)}) \geq \mathbf{y}(u_i)$, $\mathbf{y}(v_{\pi(i)}) \geq \mathbf{y}(u_j)$ and $\mathbf{y}(v_{\pi(j)}) \geq \mathbf{y}(u_j)$.\footnote{An upward planar drawing of $G$ only guarantees non-strict inequality.} Clearly, we may adjust the positions of $u_i$, $v_{\pi(i)}$, $u_j$ and $v_{\pi(j)}$ if necessary to make these inequalities strict. Thus, we may assume that $\mathbf{y}(v_{\pi(i)}) > \mathbf{y}(u_i)$, $\mathbf{y}(v_{\pi(j)}) > \mathbf{y}(u_i)$, $\mathbf{y}(v_{\pi(i)}) >\mathbf{y}(u_j)$ and $\mathbf{y}(v_{\pi(j)}) > \mathbf{y}(u_j)$.

We construct an augmented graph $\widetilde{G}_z$ from $\widetilde{G}$ as follows. Consider a point $z$ in the plane so that
\begin{equation}\label{eq.conditionz}
  \max\{\mathbf{y}(u_i),\mathbf{y}(u_j)\} < \mathbf{y}(z) < \min\{\mathbf{y}(v_{\pi(i)}), \mathbf{y}(v_{\pi(j)})\}
\end{equation}
(the strict inequalities in the previous paragraph guarantee that such a $z$ exists).

Next, as illustrated in Figure~\ref{fig.extgraph-z}, we draw four lines (shown in red) connecting $u_i$ to $z$, $u_j$ to $z$, $z$ to $v_{\pi(i)}$ and $z$ to $v_{\pi(j)}$. We add new vertices (shown in blue) to where these line-segments intersect with edges (shown in green) of $\widetilde{G}$; these vertices, together with $z$ and the vertices of $\widetilde{G}$ form the vertex set of $\widetilde{G}_z$.

The set of edges of $\widetilde{G}_z$ is obtained as follows. Include all edges of $\widetilde{G}$ that do not intersect any of the four red line-segments incident to $z$. Regard each edge of $\widetilde{G}$ that crosses the red line-segments as being subdivided by the crossing points. We also include the edges of the red line-segments which are subdivided by these crossing points.

We denote the resulting graph by $\widetilde{G}_z$. By construction, the graph $\widetilde{G}_z$ is planar and that all its edges are up-pointing.
\begin{figure}[hbt!]
    \centering
    \subfigure[]
    {\label{fig.extgraph}\includegraphics[height=0.57\textwidth]{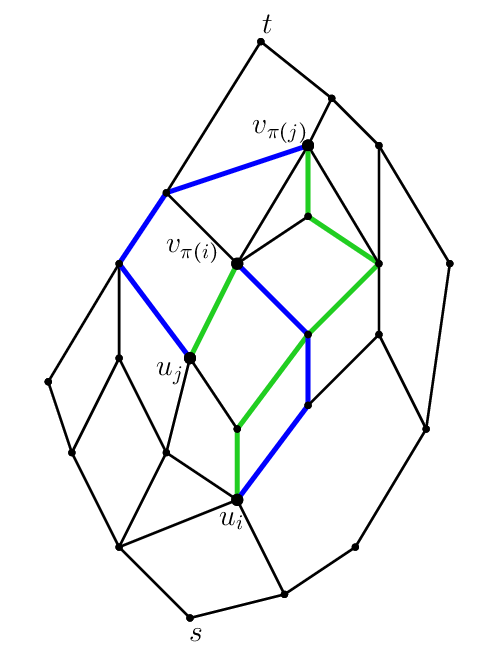}}
    \hspace{5mm}
    \subfigure[]
    {\label{fig.extgraph-z}\includegraphics[height=0.57\textwidth]{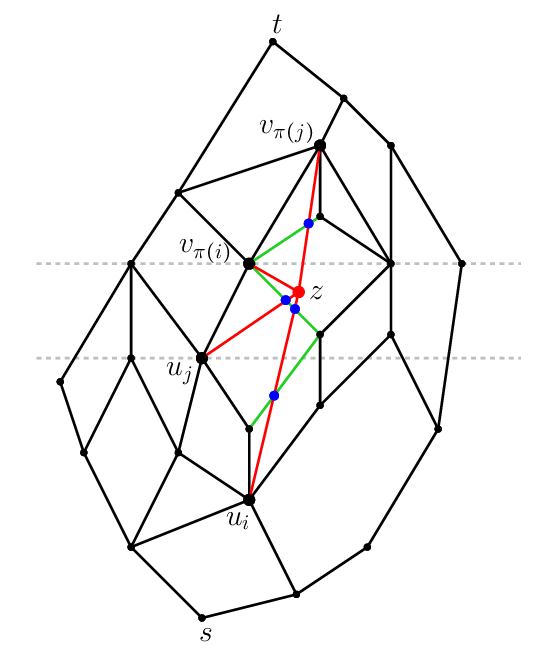}}
    \caption{(a) An illustration of two families of non-intersecting paths (in different colors) with connection types that differ by a transposition on the graph $\widetilde{G}$. (b) The augmented graph $\widetilde{G}_z$.}
    \label{fig.extendedgraph}
\end{figure}

In the augmented graph $\widetilde{G}_z$, the edges corresponding to the four red line-segments incident to $z$ naturally give rise to four new paths; denote them by
\begin{itemize}
  \item $r_i$, the path proceeding from $u_i$ to $z$ and then continuing from $z$ to $v_{\pi(i)}$
  \item $r_j$, the path from $u_j$ to $z$ and from $z$ to $v_{\pi(j)}$
  \item $r^{*}_i$, the path from $u_i$ to $z$ and from $z$ to $v_{\pi(j)}$
  \item $r^{*}_j$, the path from $u_j$ to $z$ and from $z$ to $v_{\pi(i)}$
\end{itemize}
Define two families of paths $R$ and $R^{*}$ by replacing the paths $p_i$ and $p_j$ in $P$ with the above paths as follows
\begin{equation}\label{eq.new1}
  \left\{
  \begin{array}{l}
    R = (p_1,\dotsc,p_{i-1},r_i,p_{i+1},\dotsc,p_{j-1},r_j,p_{j+1},\dotsc,p_n) \in \mathscr{P}^{\pi}(U,V) \\
    R^{*} = (p_1,\dotsc,p_{i-1},r_i^{*},p_{i+1},\dotsc,p_{j-1},r_j^{*},p_{j+1},\dotsc,p_n) \in \mathscr{P}^{\sigma\pi}(U,V)
  \end{array}
  \right.
\end{equation}

We note that $z$ is a transversal intersection (see Definition \ref{def.intersection}) of $r_i$ and $r_j$ if and only if $z$ is not a transversal intersection of $r^{*}_i$ and $r^{*}_j$. In terms of the transversal intersection numbers on the graph $\widetilde{G}_z$ (which we denote by $I_{\widetilde{G}_z}$), this implies
\begin{equation}\label{eq.swap1}
  I_{\widetilde{G}_z}(r_i^{*},r_j^{*}) = I_{\widetilde{G}_z}(r_i,r_j) \pm 1,
\end{equation}
and thus
\begin{equation}\label{eq.swap2}
  I_{\widetilde{G}_z}(R^{*}) \equiv I_{\widetilde{G}_z}(R)+1 \pmod{2}.
\end{equation}

We also note that $R^{*}$ is the image of $R$ under the involution $\phi$ described in Section~$3$. If we write $\sgn_{\widetilde{G}_z}(P)$ for the path sign of the family of paths $P$ on the graph $\widetilde{G}_z$, then by Lemma~\ref{keylemma1},  we have
\begin{equation}\label{eq.samesign}
  \sgn_{\widetilde{G}_z}(R) = \sgn_{\widetilde{G}_z}(R^{*}).
\end{equation}
Since $R$ and $P$ have the same connection type $\pi$, and $R^{*}$ and $Q$ have the same connection type $\sigma\pi$, Lemma~\ref{lem.sign+nonint} gives the following identities
\begin{equation}\label{eq.identity1}
\left\{
  \begin{array}{l}
  \sgn_{\widetilde{G}_z}(R) = \sgn_{\widetilde{G}_z}(P) (-1)^{I_{\widetilde{G}_z}(R)}  \\
  \sgn_{\widetilde{G}_z}(R^{*}) = \sgn_{\widetilde{G}_z}(Q) (-1)^{I_{\widetilde{G}_z}(R^{*})}
  \end{array}
\right.
\end{equation}

Using the results \eqref{eq.swap2} and \eqref{eq.samesign}, the identities \eqref{eq.identity1} imply
\begin{equation}\label{eq.identity3}
  \sgn_{\widetilde{G}_z}(P) = -\sgn_{\widetilde{G}_z}(Q).
\end{equation}
Then our proof will be complete once we prove the following.

\medskip
\noindent {\bf Claim.}
  Let $P,Q \in \mathscr{P}(U,V)$ be two families of paths in $\widetilde{G}$ possibly with different connection types, and not necessarily non-intersecting. Note that $P$ and $Q$ may also be viewed as families of paths in $\widetilde{G}_z$ (by replacing an edge with a subdivision of that edge, if necessary).
Then we have
  \begin{equation}\label{eq.relsign}
    \frac{\sgn_{\widetilde{G}}(P)}{\sgn_{\widetilde{G}}(Q)} = \frac{\sgn_{\widetilde{G}_z}(P)}{\sgn_{\widetilde{G}_z}(Q)}.
  \end{equation}
\medskip
\noindent {\it Proof of Claim.}
  It suffices to show that when regarding $P=(p_1,\dotsc,p_n)$ as a family of paths in $\widetilde{G}$ as opposed to a family of paths in $\widetilde{G}_z$, the path sign of $P$ changes by a factor which is independent of $P$. We remind the reader that $\sgn(P) = \prod_{i}\sgn(p_i)$ and $\sgn(p_i)$ records the parity of $|L(p_i)|$, where $L(p_i)$ is the set of marked points on the left side of $p_i$ (which recall is the region to the left of the extension of $p_i$ downward via the leftmost path originating from $s$ and upward via the leftmost path leading to $t$; see Definition~\ref{def.extension}).

\begin{figure}[htb!]
    \centering
    \subfigure[]
    {\label{fig.relsign1}\includegraphics[height=0.53\textwidth]{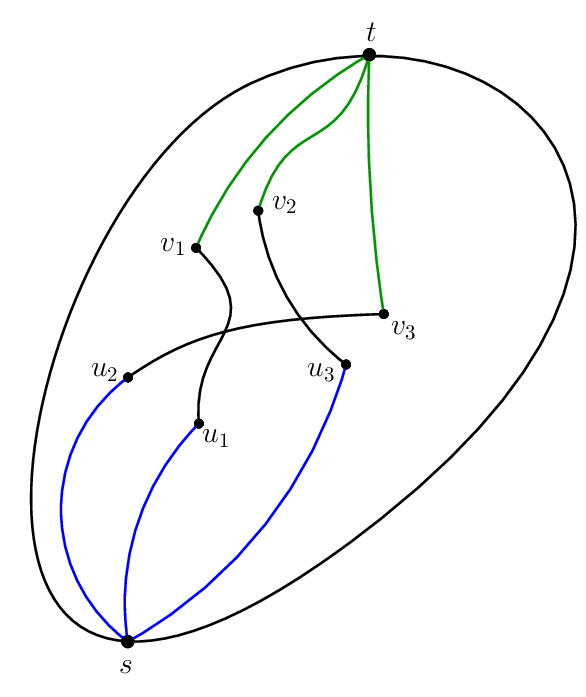}}
    \hspace{5mm}
    \subfigure[]
    {\label{fig.relsign2}\includegraphics[height=0.53\textwidth]{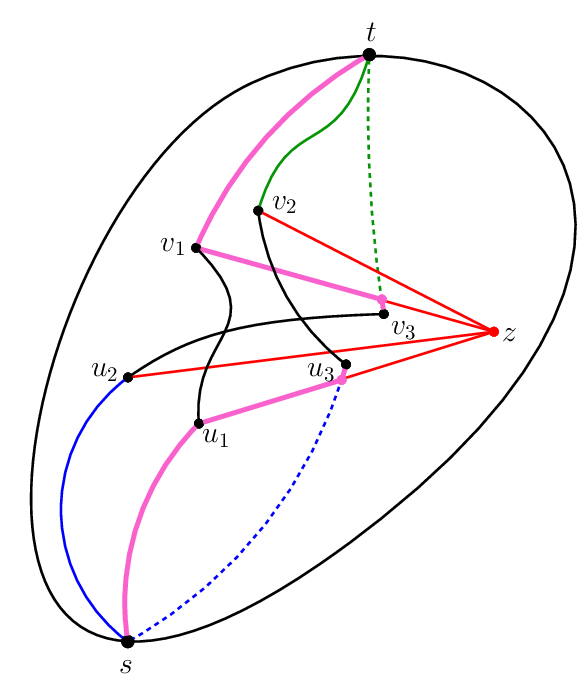}}
    \caption{(a) An example of a family of paths with $U = \{ u_1,u_2,u_3 \}$ and $V = \{ v_1,v_2,v_3 \}$. The leftmost paths from $s$ to $u_i \in U$ and from $v_i \in V$ to $t$ are shown in blue and green, respectively. (b) An illustration of how the leftmost paths from $s$ to $u_3$ and from $v_3$ to $t$ change (shown in pink) on the augmented graph $\widetilde{G}_z$.}
    \label{fig.relsign}
\end{figure}

  Note that none of the new vertices of $\widetilde{G}_z$ (compared to $\widetilde{G}$) is a marked point. However, the four red line-segments incident to $z$ in the augmented graph $\widetilde{G}_z$ may change the leftmost paths from $s$ to a starting point $u \in U$ and from an ending point $v \in V$ to $t$.

  Figure~\ref{fig.relsign2} illustrates how this may happen: the leftmost path from $s$ to $u_3$ changes, and so does the leftmost path from $v_3$ to $t$. This changes the left sides of all the paths starting from $u_3$, and also the left sides of all the paths ending at $v_3$. To be precise, the left side of each path starting at $u_3$ changes by the lower region between the dotted curve and the solid pink curve. Similarly, the left side of each path ending at $v_3$ changes by the upper region between the dotted curve and the solid pink curve.

  In general, when going from $\widetilde{G}$ to its augmented graph $\widetilde{G}_z$, there is a such region\footnote{The region reduces to the original leftmost path if there is no change on it.} $B_i$ below each starting point $u_i$, and similarly a region $B^{\prime}_i$ above each ending point $v_i$.

  The only other change that can happen in the left side of a path when regarded in $\widetilde{G}_z$ as opposed to $\widetilde{G}$ is along its left boundary. This happens when the vertex $z$ lies to the left of the whole graph (see Figure~\ref{fig.relsign_l}). In this case, the four red line-segments incident to $z$ change the leftmost path going from $s$ to $t$ by moving a portion of it further to the left. However, there are no marked points in the added region, so this has no effect on the path signs.

\begin{figure}[hbt!]
  \centering
  \includegraphics[height=0.53\textwidth]{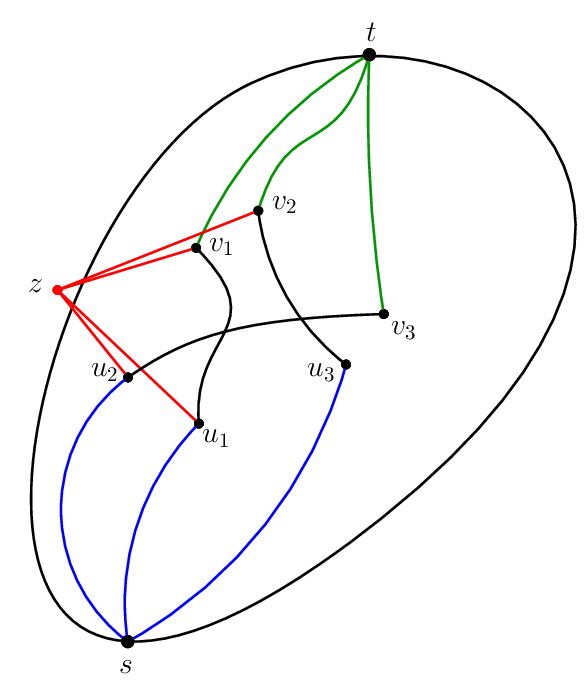}
  \caption{An illustration of the vertex $z$ lies to the left of the whole graph.}
  \label{fig.relsign_l}
\end{figure}

  The reason why our claim works is that the regions $B_i$ and $B^{\prime}_i$ are independent of the family of paths $P$. The sign of $P$ is formed by multiplying all the signs of the paths $p_i$ together. Thus, when going from $\widetilde{G}$ to $\widetilde{G}_z$, the changes in the left sides of the paths are reflected by considering all the regions $B_i$ and $B^{\prime}_i$ --- independently of $P$. This proved that the change in the path sign of $P$ when going from $\widetilde{G}$ to $\widetilde{G}_z$ is independent of $P$.
\end{proof}

\section{Counting families of non-intersecting paths connecting four collinear points}

\begin{figure}[hbt!]
  \centering
  \includegraphics[width=0.7\textwidth]{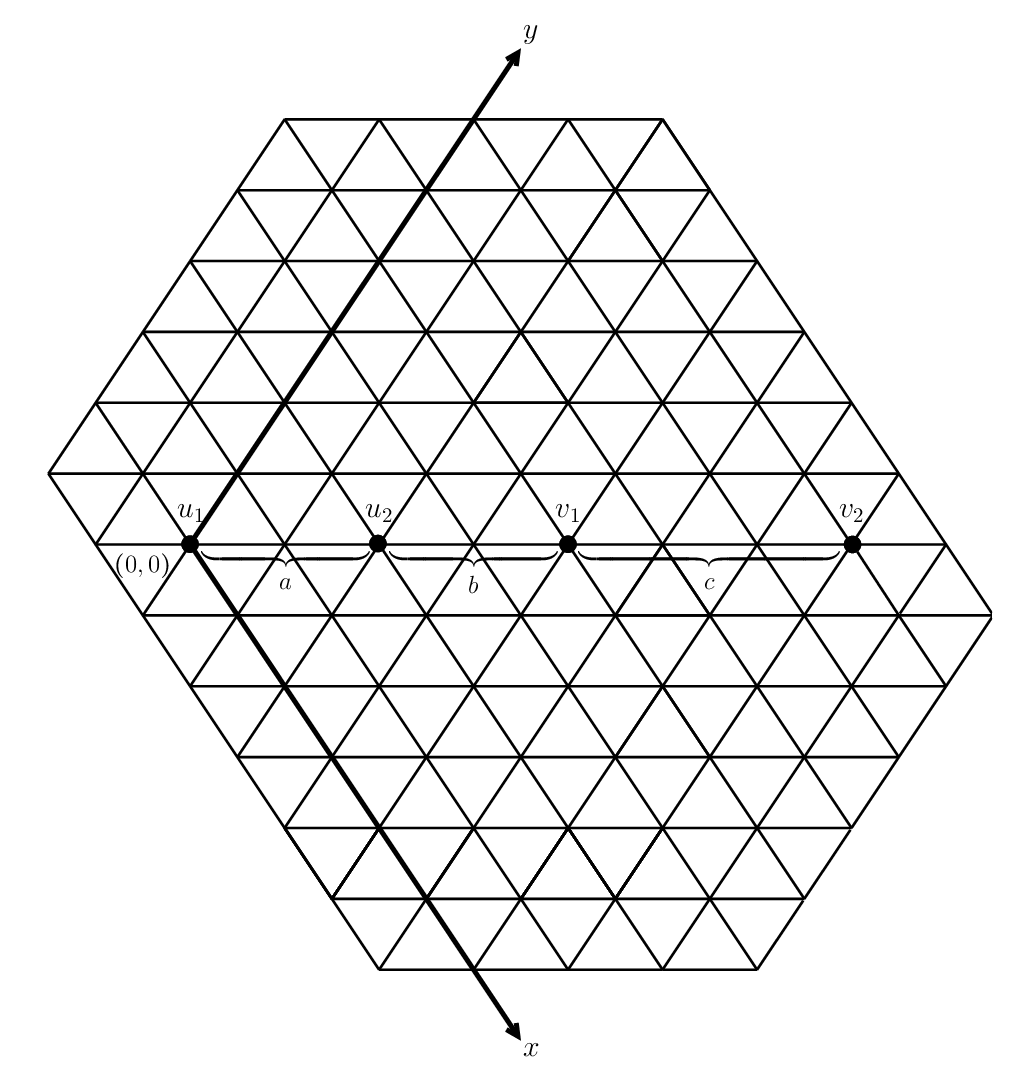}
  \caption{The coordinate system of the triangular lattice with four marked points $u_1$, $u_2$, $v_1$ and $v_2$ placed on the common horizontal line.}
  \label{fig.tri-cor}
\end{figure}

We consider an infinite triangular lattice with the orientations of each edge given in Figure~\ref{Fig1b}. Consider the weight function $\wt$ which assigns to each horizontal edge weight $z$, to each northeast pointing edge weight $y$, and to each southeast pointing edge weight $x$, where $x,y,z$ belong to some commutative ring. Choose a coordinate system on the triangular lattice (see Figure~\ref{fig.tri-cor}) by fixing a lattice point as the origin, and letting the positive $x$-axis be a lattice line pointing southeast, and the positive $y$-axis a lattice line pointing northeast.

A \textit{Delannoy path} is a lattice path from $(0,0)$ to $(n,k)$, using steps $(1,0)$, $(0,1)$ or $(1,1)$ in our coordinate system. We write $\mathscr{D}_{n,k}$ for the set of Delannoy paths going from $(0,0)$ to $(n,k)$. The weight of a path is the product of the weights of the edges in it. Let $\wt(\mathscr{D}_{n,k})$ be the total weight of each Delannoy path in $\mathscr{D}_{n,k}$. Clearly, if we assign the weights $x=y=z=1$, then $\wt(\mathscr{D}_{n,k})$ reduces to the Delannoy number, denoted by $d_{n,k}$.

The weighted Delannoy number can be defined recursively by
\begin{equation}\label{eq.Delrec}
  \wt(\mathscr{D}_{n,k}) = z \wt(\mathscr{D}_{n-1,k-1}) + x \wt(\mathscr{D}_{n-1,k}) + y \wt(\mathscr{D}_{n,k-1}),
\end{equation}
with the initial values $\wt(\mathscr{D}_{n,0}) = x^n$ and $\wt(\mathscr{D}_{0,k})=y^k$, for $n,k \geq 0$. The closed-form expression for the weighted Delannoy number is given by (see for example~\cite{FR71})
\begin{equation}\label{eq.Deleq}
  \wt(\mathscr{D}_{n,k}) = \sum_{i=0}^{n} \binom{n}{i} \binom{k}{i}x^{n-i}y^{k-i}(xy+z)^{i}.
\end{equation}

The \textit{large Schr{\"o}der path} of length $n$ is a lattice path going from $(0,0)$ to $(n,n)$ that never passes below the line $y=x$, consisting of the steps $(1,0)$, $(0,1)$ or $(1,1)$ in our coordinate system. Let $\mathscr{R}_n$ be the set of large Schr{\"o}der paths of length $n$ and $\wt(\mathscr{R}_n)$ the total weight of paths in $\mathscr{R}_n$. It is obvious that the large Schr{\"o}der number $r_n$ is obtained from $\wt(\mathscr{R}_n)$ by assigning the weights $x=y=z=1$. The following simple formula is given in \cite{CP17}:
\begin{equation}\label{eq.Sch}
  \wt(\mathscr{R}_n) = \sum_{i=0}^{n}\frac{1}{i+1}\binom{2i}{i}\binom{n+i}{n-i}x^iy^iz^{n-i}.
\end{equation}

We study the special case when $U=\{u_1,u_2\}$, $V=\{v_1,v_2\}$ and these four points are placed on the line $y=x$, at relative spacings given by three parameters $a,b$ and $c$, where $a,c \in \mathbb{Z}_{>0}$ and $b \in \mathbb{Z}_{\geq 0}$ (see Figure~\ref{fig.tri-cor}). The point $u_1$ is fixed at $(0,0)$, $u_2$ is at $(a,a)$, $v_1$ is at $(a+b,a+b)$ and $v_2$ is at $(a+b+c,a+b+c)$. Let $\mathscr{N}_{0}(a,b,c)$ be the set of non-intersecting paths from $U$ to $V$ on the triangular lattice. Then we have the following result.
\begin{theorem}\label{thm.2patht}
  The total weight of the non-intersecting paths described above is
  \begin{equation*}
    \wt(\mathscr{N}_0(a,b,c)) = 2 \wt(\mathscr{D}_{b,b}) \left( \sum_{i=1}^{a} \sum_{j=1}^{c} \wt(\mathscr{R}_{b+i+j-1}) \wt(\mathscr{D}_{a-i,a-i}) \wt(\mathscr{D}_{c-j,c-j})) \right),
  \end{equation*}
where the weights involved on the right hand side have the explicit expressions given by formulas~\eqref{eq.Deleq} and~\eqref{eq.Sch}.
\end{theorem}
\begin{proof}
  Let $G$ be the subgraph of the oriented triangular lattice induced by the vertices in the rhombus whose four vertices are $(0,0),(a+b+c,0),(a+b+c,a+b+c)$ and $(0,a+b+c)$ in our coordinate system. We will apply Theorem~\ref{thm.main}, with the $st$-planar graph $\widetilde{G}$ chosen to be $G$, with the source $u_1$ and the sink $v_2$. An upward planar drawing of $G$ is obtained by rotating this rhombus counterclockwise $90$ degrees.

  We will in fact find it more convenient not to perform this rotation, but instead think of Theorem~\ref{thm.main} as phrased for graphs having a ``rightward'' planar drawing; the left side of a path becomes then the ``top'' side from this point of view.

  We note that the set of starting points $\{u_1,u_2\}$ and the set of ending points $\{v_1,v_2\}$ are not compatible (in the sense of the paragraph after Theorem~\ref{thm.GV1}) in our situation, so the Lindstr{\"o}m-Gessel-Viennot theorem does not apply. However, by our main result, the enumeration of such families of non-intersecting paths is given by the determinant of the matrix $M$ in Theorem~\ref{thm.main}. We compute the entries of the matrix $M$ as follows.

  The $(2,1)$-entry: The paths from $u_2$ to $v_1$. All the four marked points are contained in the top side of each path. The total weight of these paths is $\wt(\mathscr{D}_{b,b})$, so the $(2,1)$-entry of the matrix $M$ is $(-1)^{4}\wt(\mathscr{D}_{b,b})$.

  The $(1,1)$-entry: The paths from $u_1$ to $v_1$. We first consider the set of paths that do not pass through the point $u_2$; denote this by $\mathscr{P}(u_1,v_1;\widehat{u_2})$. Define an involution $\gamma_{u_2}$ on this set, called the \textit{involution with respect to $u_2$}, as follows.
  \begin{enumerate}
    \item Consider a path $p$ of this set. Find the nearest vertex to $u_2$ in which the path $p$ intersects the line $y=x$ to the left (resp., right) of the point $u_2$, say $\alpha$ (resp., $\beta$).
    \item Reflect the subpath of $p$ from $\alpha$ to $\beta$ across the line $y=x$ and keep the other parts of $p$ fixed.
  \end{enumerate}
  For each path $p \in \mathscr{P}(u_1,v_1;\widehat{u_2})$, this yields a new path $\gamma_{u_2}(p)\in\mathscr{P}(u_1,v_1;\widehat{u_2})$ with the opposite path sign, because the sets $L(p)$ and $L(\gamma_{u_2}(p))$ differ just in that one of them contains the marked point $u_2$, while the other does not. Since the subpath of the path $p$ from $\alpha$ to $\beta$ consists of the same number of $(1,0)$ steps and $(0,1)$ steps, this involution is weight-preserving. This implies that in the $(1,1)$-entry of the matrix $M$, the terms corresponding to the paths in $\mathscr{P}(u_1,v_1;\widehat{u_2})$ are cancelled out.

  On the other hand, the total weight of the paths that contain the point $u_2$ is given by $\wt(\mathscr{D}_{a,a}) \wt(\mathscr{D}_{b,b})$, and their path signs are $(-1)^{4}$. Thus, the $(1,1)$-entry of the matrix $M$ is $(-1)^{4} \wt(\mathscr{D}_{a,a}) \wt(\mathscr{D}_{b,b})$.

  The $(2,2)$-entry: The paths from $u_2$ to $v_2$. This is similar to Case $2$ (apply the involution $\gamma_{v_1}$ on the set $\mathscr{P}(u_2,v_2;\widehat{v_1})$). The $(2,2)$-entry of the matrix $M$ is then $(-1)^{4}\wt(\mathscr{D}_{b,b})\wt(\mathscr{D}_{c,c})$.

  The $(1,2)$-entry: The paths from $u_1$ to $v_2$. It is not hard to verify (using the involutions with respect to $u_2$ and $v_1$ described above) that in the expression of Theorem~\ref{thm.main} for the $(1,2)$-entry of the matrix $M$, all terms cancel out except the terms corresponding to the following two sets of paths connecting $u_1$ to $v_2$:

  \begin{enumerate}
    \item The set of paths that pass through both $u_2$ and $v_1$. The path sign of each of these paths is $(-1)^{4}$. The total weight of this set is therefore $\wt(\mathscr{D}_{a,a})\wt(\mathscr{D}_{b,b})\wt(\mathscr{D}_{c,c})$.
    \item The set of the paths that do not touch the line $y=x$ in between the marked points $u_2$ and $v_1$. The paths that pass above both $u_2$ and $v_1$ can be viewed as the concatenation of five (possibly empty) subpaths:
        \begin{itemize}
          \item from $u_1$ to $(a-i,a-i)$ via a Delannoy path,
          \item the unit step from $(a-i,a-i)$ to $(a-i,a-i+1)$,
          \item from $(a-i,a-i+1)$ to $(a+b+j-1,a+b+j)$ via a Schr{\"o}der path of length $b+i+j-1$.
          \item the unit step from $(a+b+j-1,a+b+j)$ to $(a+b+j,a+b+j)$, and
          \item from $(a+b+j,a+b+j)$ to $v_2$ via a Delannoy path,
        \end{itemize}
        where $1 \leq i \leq a$ and $1 \leq j \leq c$.

        Writing $\rho$ for the total weight of the paths that pass above the points $u_2$ and $v_1$, this implies that $\rho$ is given by
        \begin{equation}\label{eq.M1}
          \rho = \sum_{i=1}^{a}\sum_{j=1}^{c} \wt(\mathscr{R}_{b+i+j-1}) \wt(\mathscr{D}_{a-i,a-i}) \wt(\mathscr{D}_{c-j,c-j}).
        \end{equation}
        Each of these paths has path sign $(-1)^{4}$.

        By symmetry, the paths that pass below both $u_2$ and $v_1$ have the same total weight $\rho$. Since their path signs are $(-1)^{2}$, the total weight of the paths in this subcase is then $2\rho$.
  \end{enumerate}
  Therefore, the $(1,2)$-entry of the matrix $M$ is $\wt(\mathscr{D}_{a,a})\wt(\mathscr{D}_{b,b})\wt(\mathscr{D}_{c,c}) + 2\rho$.

  By Theorem~\ref{thm.main}, the total weight of the families of non-intersecting paths connecting $\{u_1,u_2\}$ to $\{v_1,v_2\}$ is then
  \begin{equation*}
  \left| \det M \right| = \left| \det
    \begin{pmatrix}
      \wt(\mathscr{D}_{a,a})\wt(\mathscr{D}_{b,b}) & \wt(\mathscr{D}_{a,a})\wt(\mathscr{D}_{b,b})\wt(\mathscr{D}_{c,c})+ 2\rho \\
      \wt(\mathscr{D}_{b,b}) & \wt(\mathscr{D}_{b,b})\wt(\mathscr{D}_{c,c})
    \end{pmatrix}
    \right|
   = 2\wt(\mathscr{D}_{b,b})\rho.
  \end{equation*}
This completes the proof of Theorem~\ref{thm.2patht}.
\end{proof}

We can readily deduce a counterpart result on the square lattice by setting the weight $z=0$ in Theorem~\ref{thm.2patht}. The weighted Delannoy number is replaced by the weighted binomial coefficient while the weighted large Schr{\"o}der number by the weighted Catalan number. This leads to the following result.
\begin{corollary}\label{cor.2paths}
  Let $u_1=(0,0)$, $u_2=(a,a)$, $v_1=(a+b,a+b)$ and $v_2=(a+b+c,a+b+c)$ be four points on the square lattice $\mathbb{Z}^2$, where $a,c\in\mathbb{Z}_{>0}$ and $b\in\mathbb{Z}_{\geq 0}$.

  The total weight of non-intersecting paths from $u_k$ to $v_{\ell}$, $k,\ell=1,2$, on the square lattice is
  \begin{equation*}
    \wt(\mathscr{N}_{0}(a,b,c)) |_{z=0} = 2 \wt(\mathscr{B}_{b,b}) \left( \sum_{i=1}^{a}\sum_{j=1}^{c} \wt(\mathscr{C}_{b+i+j-1}) \wt(\mathscr{B}_{a-i,a-i}) \wt(\mathscr{B}_{c-j,c-j}) \right),
  \end{equation*}
  where $\wt(\mathscr{B}_{n,k}) = \binom{n+k}{k}x^{k}y^{n}$ is the weighted binomial coefficient and $\wt(\mathscr{C}_n) = \frac{1}{n+1}\binom{2n}{n} x^ny^n$ is the weighted Catalan number.
\end{corollary}

\section{Enumeration of domino tilings of a mixed Aztec rectangle with unit holes}

\begin{figure}[hbt!]
    \centering
    \subfigure[]
    {\label{fig.AD}\includegraphics[width=0.3\textwidth]{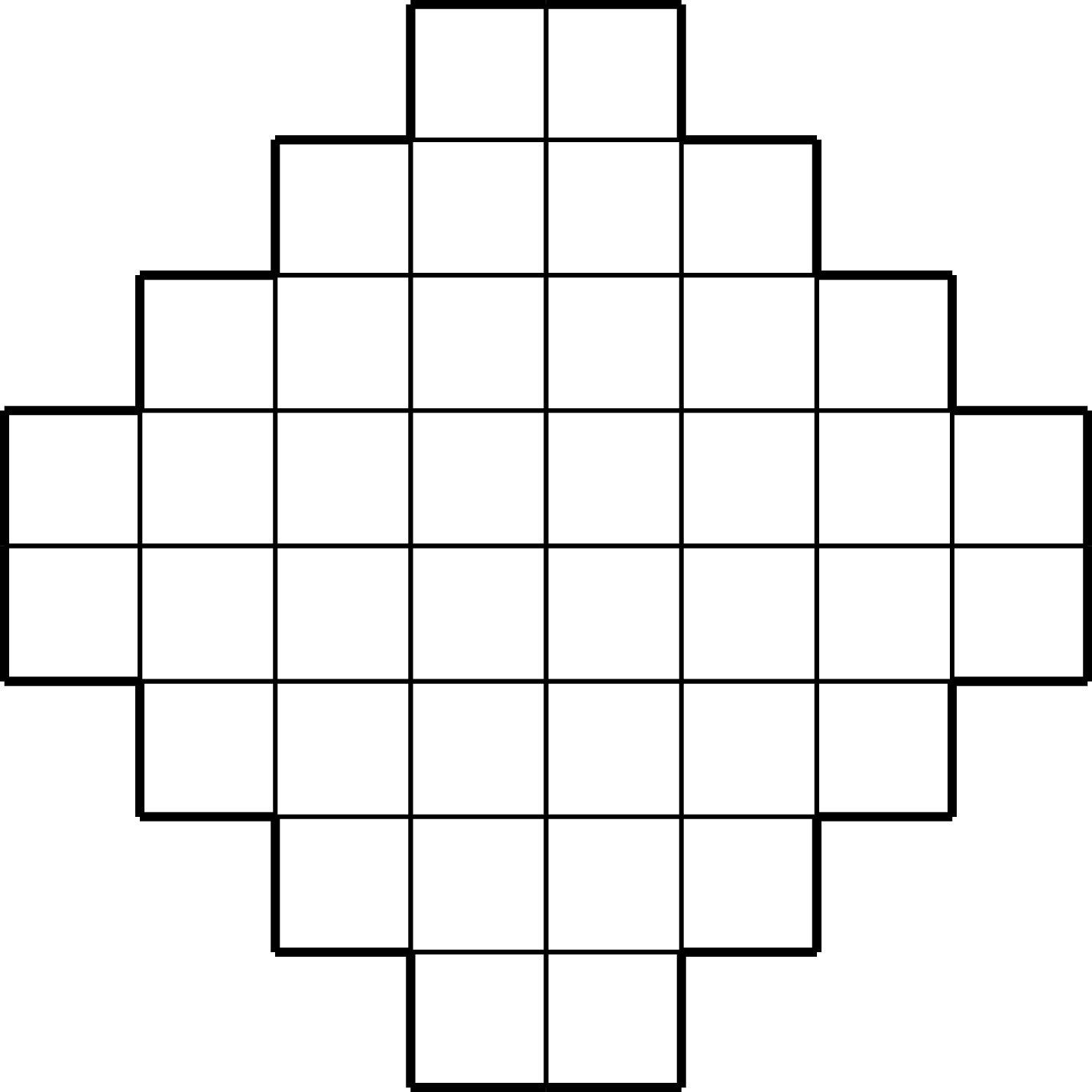}}
    \hspace{20mm}
    \subfigure[]
    {\label{fig.AR}\includegraphics[width=0.3\textwidth]{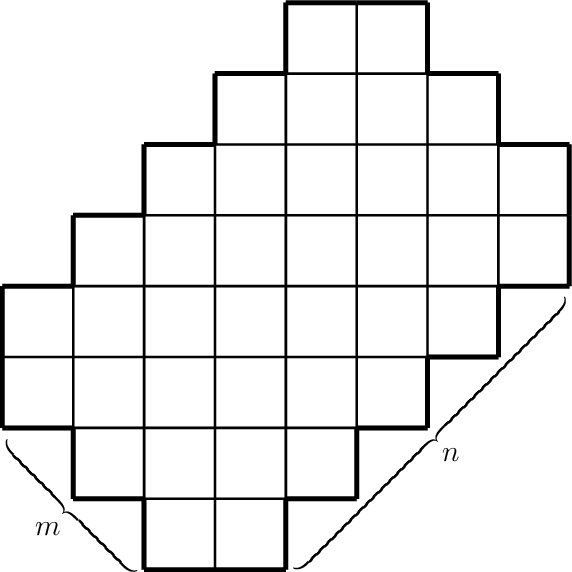}}
    \caption{(a) The Aztec diamond of order $4$. (b) The Aztec rectangle with $m=3$ and $n=5$.}
    \label{fig.ADAR}
\end{figure}

The \textit{Aztec diamond} of order $n$ (see Figure~\ref{fig.AD}) is the union of all unit squares in the region  $|x|+|y| \leq n+1$. The \textit{Aztec rectangle} $AR_{m,n}$ is a natural generalization of the Aztec diamond; there are $m$ unit squares on the southwest side and $n$ unit squares on the southeast side, as shown in Figure~\ref{fig.AR}. The Aztec diamond is the special case when $m=n$.

In the past three decades, there have been numerous results on the enumeration of domino tilings of Aztec rectangle regions with holes. In the early survey paper ~\cite[Section 3]{Propp99}, Propp listed some problems about the domino tilings with certain holes. We refer the reader to \cite{Ciucu97}, \cite{HI99}, \cite{Kra00}, \cite{Lai15} and \cite{Sai17} for the enumeration formulas and the different approaches to these domino tiling problems. In \cite{Ciucu14} and \cite{Ciucu16}, Ciucu investigated the interactions of different types of monomer clusters (sets of unit holes) on a two dimensional dimer system and determined the asymptotic of their correlations.

\begin{figure}[hbt!]
    \centering
    \subfigure[]
    {\label{fig.ADm}\includegraphics[width=0.3\textwidth]{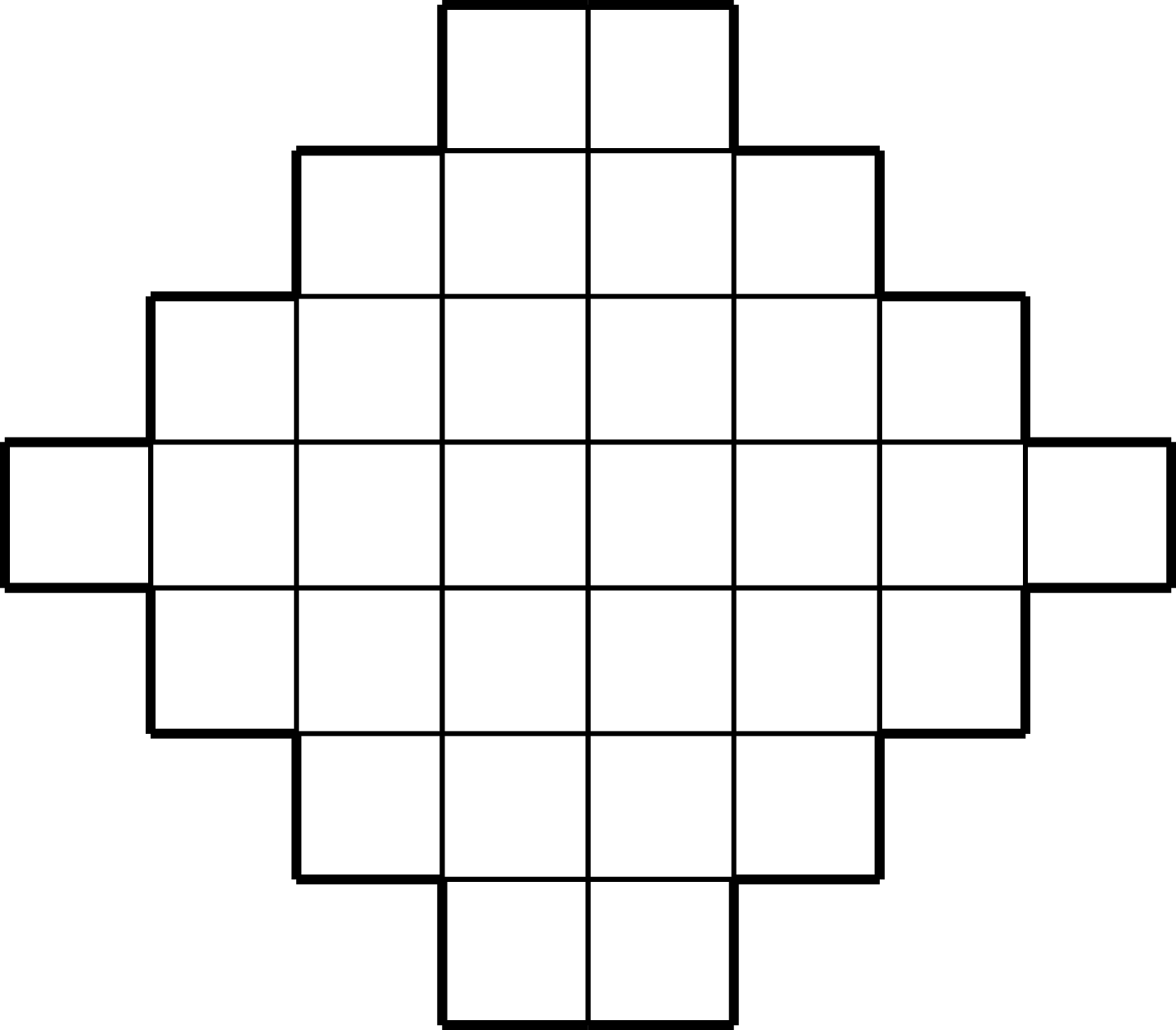}}
    \hspace{20mm}
    \subfigure[]
    {\label{fig.ARm}\includegraphics[width=0.3\textwidth]{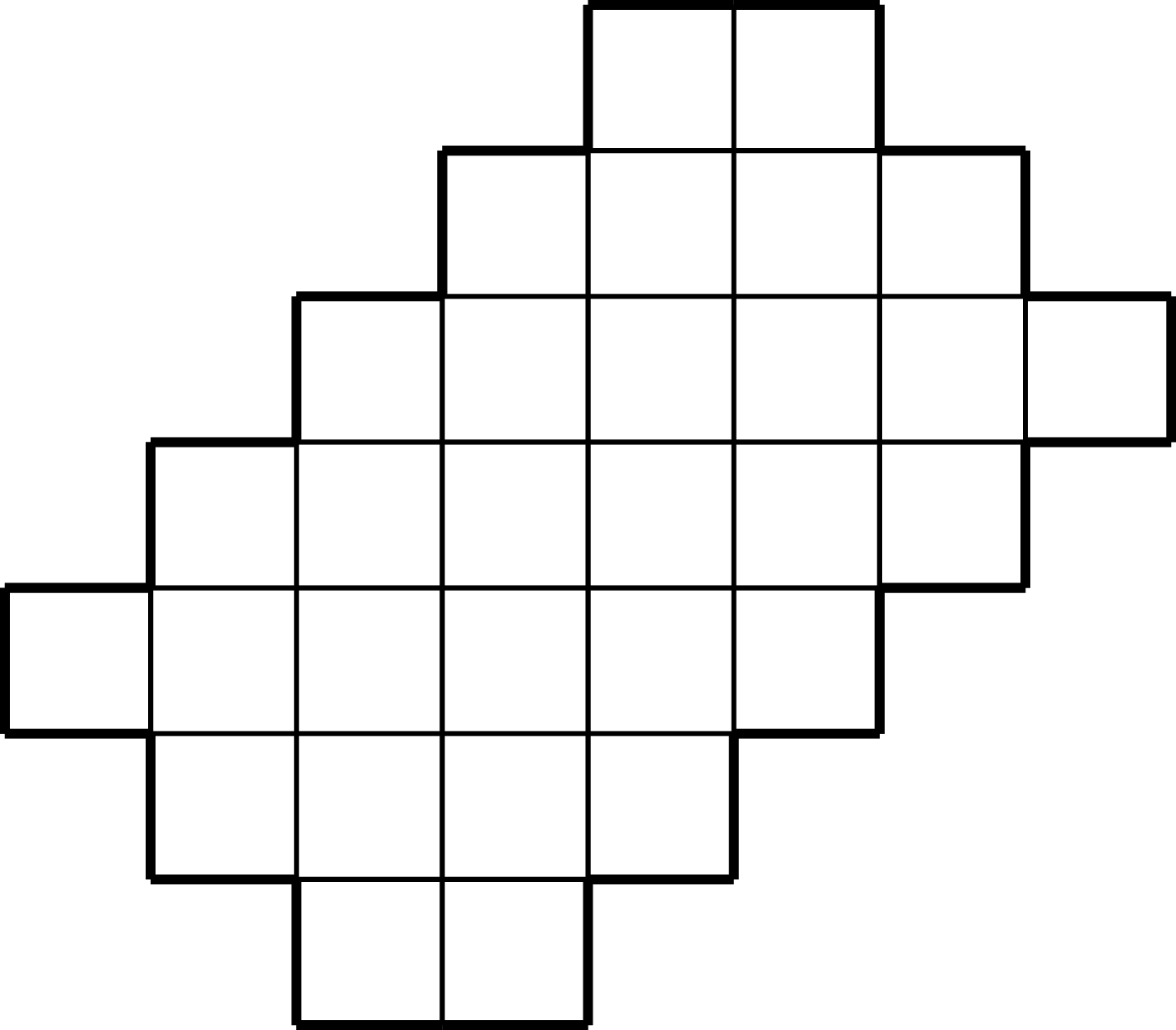}}
    \caption{(a) The mixed Aztec diamond of order $4$. (b) The mixed Aztec rectangle with $m=3$ and $n=5$.}
    \label{fig.ADmARm}
\end{figure}

In~\cite[Section 3]{Ciucu08}, Ciucu introduced the \textit{mixed Aztec diamond} of order $n$, $MD_n$, which can be viewed as the Aztec diamond of order $n$ with the southwest $n$ unit squares and the southeast $n$ unit squares removed (see Figure~\ref{fig.ADm}). Similarly, the \textit{mixed Aztec rectangle} $MR_{m,n}$ is defined to be the Aztec rectangle $AR_{m,n}$ with the southwest $m$ unit squares and the southeast $n$ unit squares deleted (see Figure~\ref{fig.ARm}).

It is easy to see that there is only one domino tiling of the mixed Aztec rectangle with any positive integers $m$ and $n$; all the dominoes are horizontal.

Given a mixed Aztec rectangle, we consider the checkerboard coloring of the square lattice with the unit squares along its top right side colored black.

\begin{theorem}\label{thm.invariant}
  The number of domino tilings of the mixed Aztec rectangle with arbitrary unit holes is invariant under color-preserving translations of the set of holes\footnote{A translation of a set of holes is obtained by applying the same translation to each hole in the set.}, provided all the unit holes are still contained in the mixed Aztec rectangle.
\end{theorem}
\begin{proof}
  In the checkerboard coloring mentioned before the statement of Theorem~\ref{thm.invariant}, mark the midpoint of the left edge of each black unit square, and join these midpoints by edges as shown in Figure~\ref{fig.mARgrid}, to obtain a subgraph of the triangular lattice.

  \begin{figure}[hbt!]
    \centering
    \subfigure[]
    {\label{fig.mARgrid}\includegraphics[height=0.38\textwidth]{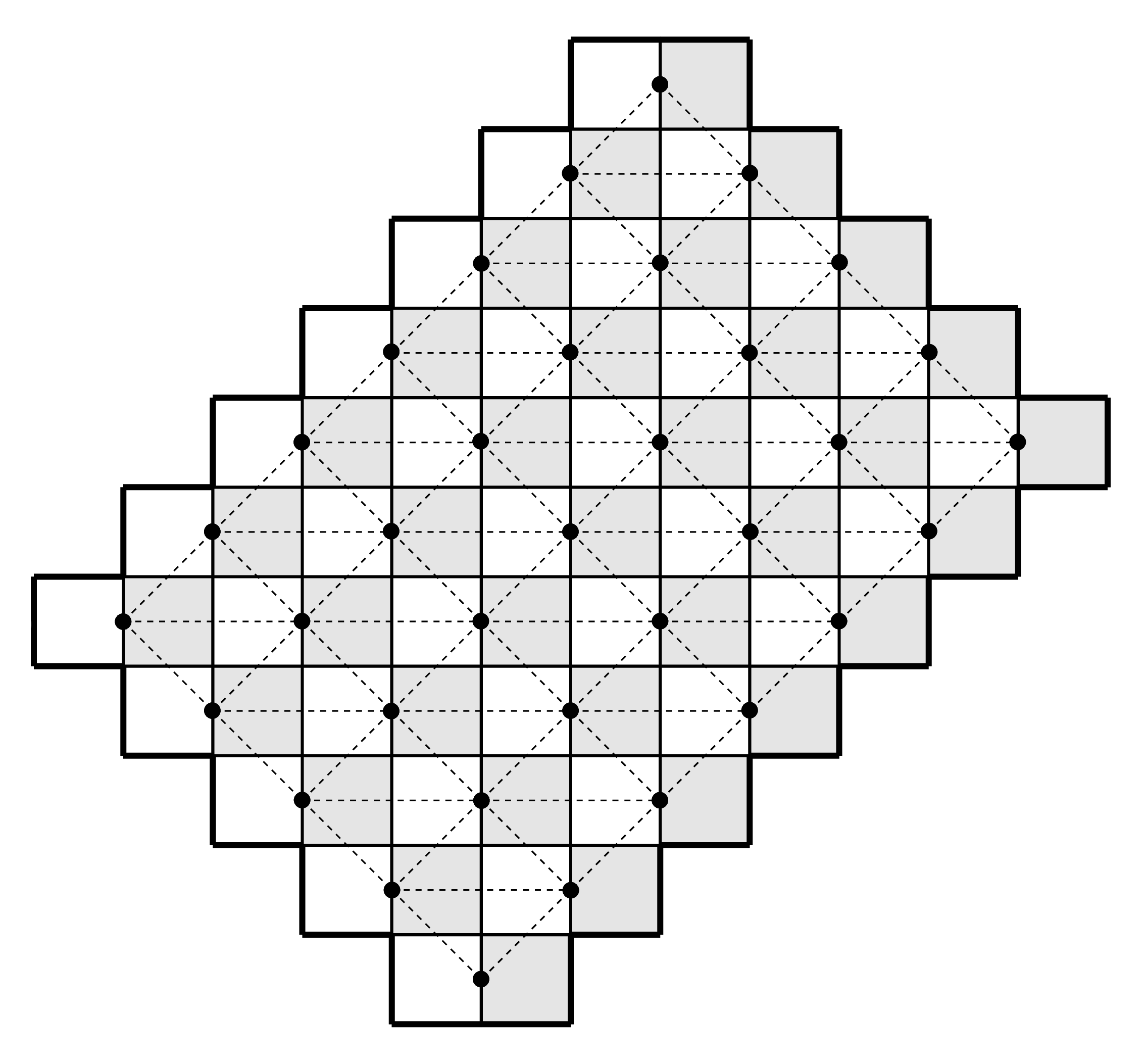}}
    \hspace{10mm}
    \subfigure[]
    {\label{fig.mARholes}\includegraphics[height=0.38\textwidth]{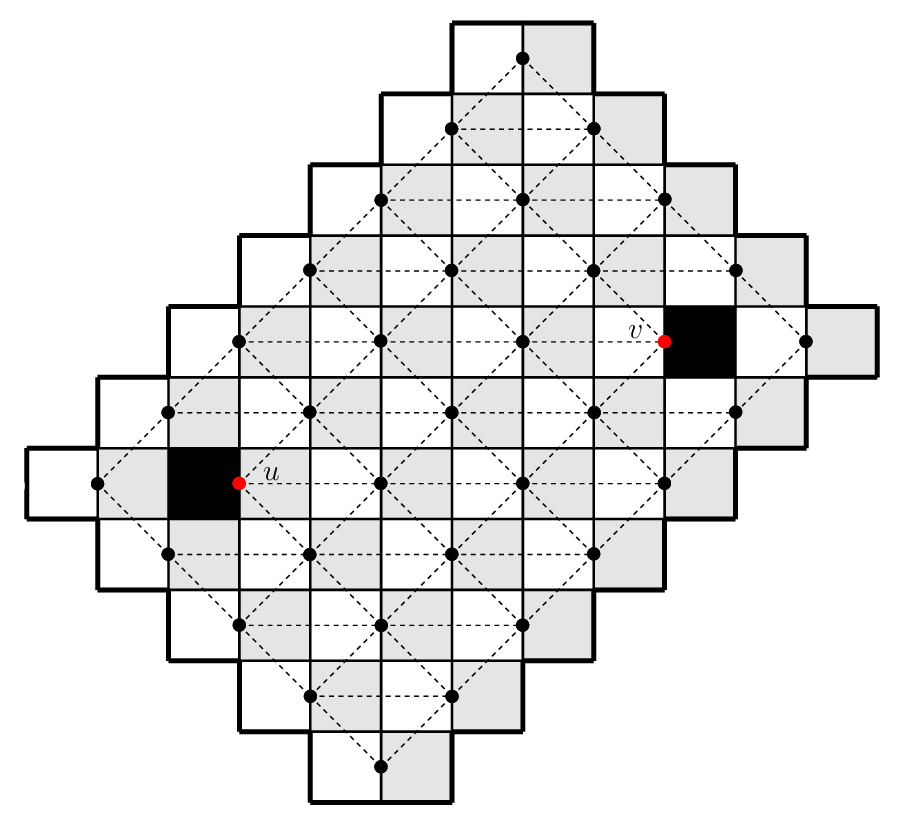}}
    \caption{(a) The checkerboard coloring of the mixed Aztec rectangle and a subgraph of the triangular lattice (dotted edges) on it. (b) An illustration of the white and black unit holes and their corresponding endpoints (shown in red) of paths. }
    \label{fig.mAR}
  \end{figure}

  We use the fact that there is a bijection between the set of domino tilings of a region $R$ on the square lattice and families of non-intersecting Delannoy paths with certain starting and ending points (determined by the region $R$). This construction is implicit in work of Sachs and Zernitz~\cite{SZ94}, and was made explicit by Randall (see \cite[Section 4]{Ciucu96} and~\cite[Section 2]{LRS01}).

  The bijection works as follows (see Figure~\ref{fig.mARex}). Given a domino tiling, map each vertical domino with the top (resp., bottom) unit square colored black to a $(1,0)$ (resp., $(0,1)$) step on the triangular lattice, and each horizontal domino with the left unit square colored black to a $(1,1)$ step on the triangular lattice. Note that no step of the lattice paths corresponds to horizontal dominos in which the right unit square is black.

  Note that each unit hole generates one starting or ending point of a lattice path (see Figure~\ref{fig.mARholes}). More precisely, under the checkerboard coloring, if a unit hole is colored white, then the midpoint on its right edge becomes a starting point $u$ of a lattice path. Similarly, if a unit hole is colored black, then the midpoint on its left edge becomes an ending point $v$ of a lattice path. Figure~\ref{fig.mARex} shows an example of a domino tiling of the mixed Aztec rectangle with two unit holes and the corresponding lattice path.

  \begin{figure}[hbt!]
    \centering
    \subfigure[]
    {\label{fig.mARex2}\includegraphics[height=0.38\textwidth]{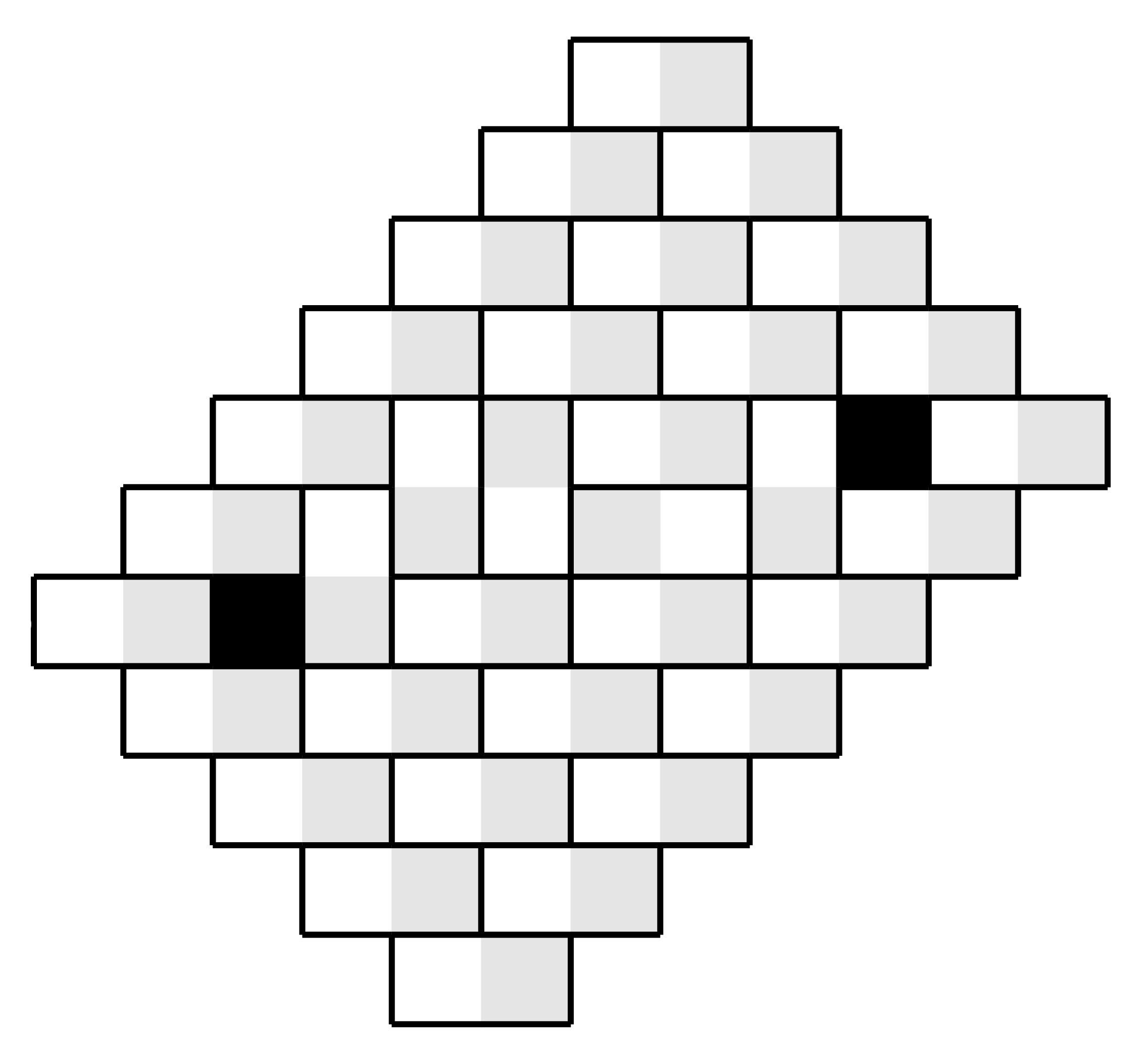}}
    \hspace{10mm}
    \subfigure[]
    {\label{fig.mARex1}\includegraphics[height=0.38\textwidth]{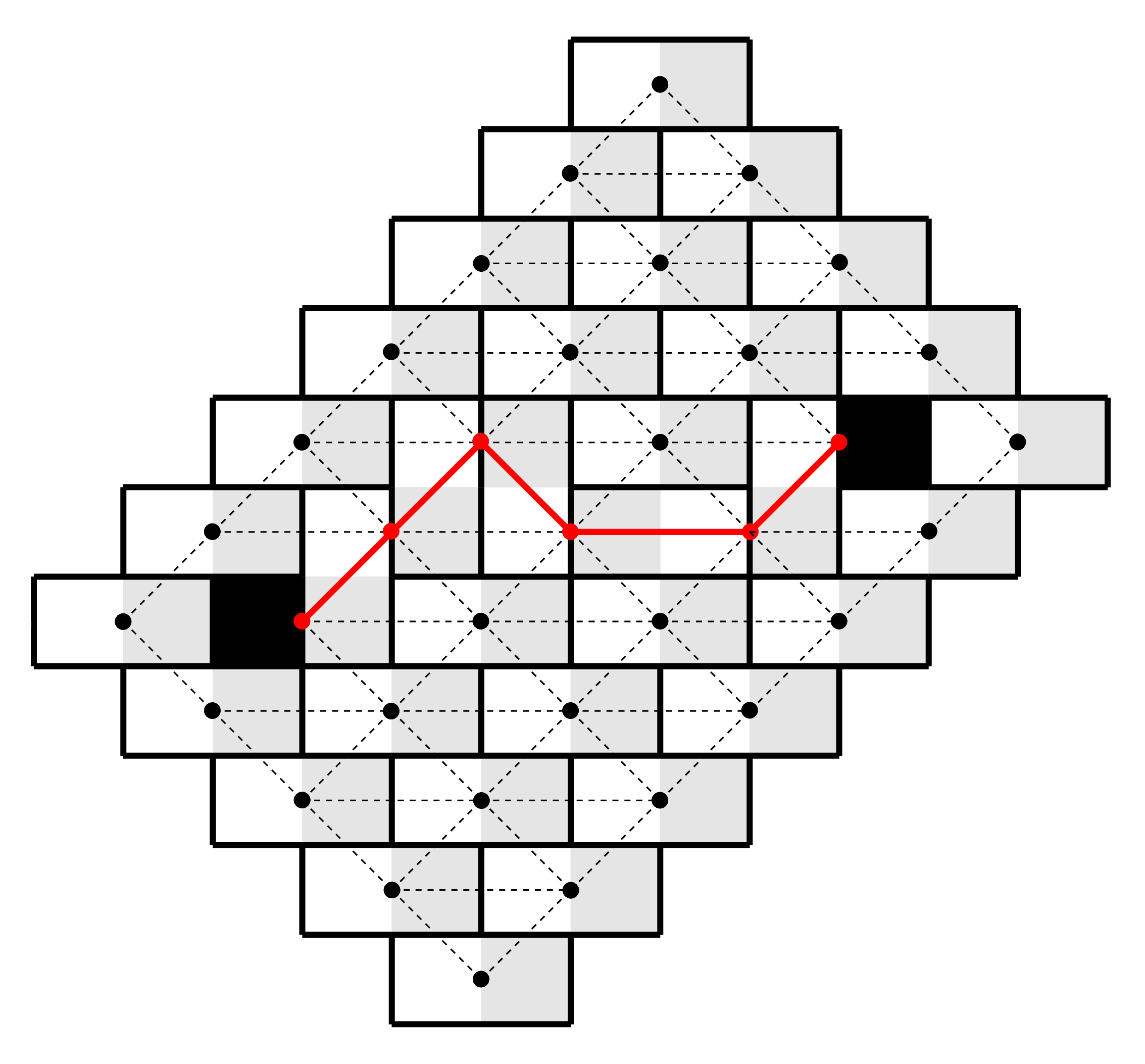}}
    \caption{(a) An example of a domino tiling of the mixed Aztec rectangle with two unit holes. (b) The corresponding lattice path from Figure~\ref{fig.mARex2}.}
    \label{fig.mARex}
  \end{figure}

  Consider now the mixed Aztec rectangle with arbitrary unit holes and the subgraph of the triangular lattice induced on it by the above bijection. Let $U=\{u_1,\dotsc,u_k\}$ be the starting points of paths corresponding to the white unit holes, and $V=\{v_1,\dotsc,v_{\ell}\}$ the set of ending points corresponding to the black unit holes. We note that there is no tiling if $k \neq \ell$.

  By the above bijection, the domino tilings of the mixed Aztec rectangle with these holes are in bijection with the families of non-intersecting Delannoy paths connecting $U$ to $V$.

  Under a color-preserving translation of the unit holes, assuming their images are still contained in the mixed Aztec rectangle, the cardinalities of $U$ and $V$ do not change. Furthermore, for each pair $(i,j)$, the relative position of the points $u_i$ and $v_j$ on the triangular lattice does not change. Thus, the resulting families of non-intersecting paths before and after the translation are translates of one another. By the bijection between domino tilings and families of non-intersecting Delannoy paths, it follows that the number of domino tilings after the translation is the same as before.
\end{proof}

It seems difficult to find formulas for the number of tilings of mixed Aztec rectangles with unit holes in general positions. We present below a simple special case in which we can give such a formula.

Consider the case when we have four unit holes along a common horizontal line inside the mixed Aztec rectangle. It follows from the above mentioned bijection that this region has no tilings unless the color sequence of these holes, from left to right, is white-black-white-black or white-white-black-black. In the former case, the family of paths corresponding to the tilings consists of two independent Delannoy paths, so the number of tilings readily follows to be the product of two Delannoy numbers. The latter case is more interesting, and we discuss it in the rest of this section.

Suppose the spacing between the holes are (from left to right) $2a-1$, $2b$ and $2c-1$ units. Let ${MR}_{m,n}(a,b,c)$ be the mixed Aztec rectangle region with holes described above (we always assume that the values of the parameters and the placement of the holes are so that the four holes are inside the mixed Aztec rectangle). Then we have the following theorem.

\begin{theorem}
  The number of domino tilings of the regions ${MR}_{m,n}(a,b,c)$ is given by
  \begin{equation}\label{eq.mAR4}
    \M({MR}_{m,n}(a,b,c)) = 2 d_{b,b} \left( \sum_{i=1}^{a} \sum_{j=1}^{c} r_{b+i+j-1} d_{a-i,a-i}d_{c-j,c-j} \right),
  \end{equation}
  where $d_{n,k}$ is the Delannoy number and $r_n$ is the large Schr{\"o}der number\footnote{The numbers $d_{n,k}$ and $r_n$ are obtained from equations~\eqref{eq.Deleq} and~\eqref{eq.Sch} by setting the weights $x=y=z=1$, respectively.}.

  In particular, $\M({MR}_{m,n}(a,b,c))$ only depends on the separations of the four holes (and not on $m$, $n$, or the position of the left hole).
\end{theorem}
%


%
\begin{proof}
    According to the bijection described in the proof of Theorem~\ref{thm.invariant}, the domino tilings of ${MR}_{m,n}(a,b,c)$ are in bijection with the families of Delannoy paths having the set of starting points $U = \{ u_1,u_2 \}$ and ending points $V = \{ v_1,v_2 \}$ on the triangular lattice described above.

    These four points lie on a common horizontal line. From left to right, they are $u_1$, $u_2$, $v_1$ and $v_2$ with $a$, $b$ and $c$ $(1,1)$-steps between each consecutive pair, respectively. We point out that this setting is the same as the one in Section $5$. Using Theorem~\ref{thm.2patht} and setting the weights $x=y=z=1$, we immediately obtain formula \eqref{eq.mAR4}. By this formula, the number of domino tilings is clearly independent of $m$, $n$ and the position of the left unit hole (this also follows directly from Theorem~\ref{thm.invariant}).

\end{proof}
%


\section{Concluding remarks}
In this paper, we presented a determinant formula for the enumeration of families of non-intersecting paths with arbitrary starting and ending points on a directed acyclic graph having an upward planar drawing (Theorem~\ref{thm.main}). In the special case when the starting and ending points satisfy a certain restrictive condition, this is given by the Lindstr{\"o}m-Gessel-Viennot theorem (even for directed acyclic graph that do not have an upward planar drawing). The challenge was to address the case of arbitrary starting and ending points. It is natural to ask whether our result can be extended to arbitrary directed acyclic graphs.

In general, the entries of matrix $M$ in~\eqref{eq.main} are complicated and $\det M$ seems difficult to evaluate; only some simple examples were worked out in Section~$5$. We hope our result can shed new light on studies concerning the Lindstr{\"o}m-Gessel-Viennot theorem, especially the enumeration of tilings of a region with holes.

\subsection*{Acknowledgements}

The author thanks Mihai Ciucu for stimulating discussions and helpful suggestions on the preliminary version of this paper. The author also thanks the reviewers for carefully reading the manuscript and providing helpful comments.



\end{document}